\documentclass[final]{siamart190516}
\usepackage{cite}
\usepackage[utf8]{inputenc}
\usepackage[T1]{fontenc}
\usepackage[english]{babel}
\usepackage{amsfonts}
\usepackage{amsmath}
\usepackage{amssymb}
\usepackage{amsopn}
\usepackage{hyperref}[6.83]
\usepackage{algorithm}
\usepackage{algorithmicx}
\usepackage[noend]{algpseudocode}
\usepackage[title]{appendix}
\usepackage{bm}
\usepackage{booktabs}
\usepackage{dsfont}
\usepackage{enumitem}
\usepackage{fancyhdr}
\usepackage{graphicx}
\usepackage{grffile}
\usepackage{import}
\usepackage{mathrsfs}
\usepackage{mathdots}
\usepackage{mathtools}
\usepackage{microtype}
\usepackage{array}
\usepackage{overpic}
\usepackage{varwidth}
\usepackage{xcolor}
\usepackage{listings}
\usepackage{import}
\usepackage[caption=false]{subfig}
\usepackage{tikz}
\usetikzlibrary{shapes,arrows,plotmarks,positioning,calc}
\usetikzlibrary{arrows.meta}
\usepackage{pgfplots}
\pgfplotsset{compat=newest}
\usepgfplotslibrary{patchplots}
\usepackage{todonotes}
\usepackage{tabularx}
\usepackage{xspace}
\usepackage{matlab-prettifier}

\addto\captionsenglish{\renewcommand{\figurename}{Fig.}}

\makeatletter
\algrenewcommand\ALG@beginalgorithmic{\small}
\makeatother

\def\intcolor{cyan!60}
\def\extcolor{red!60}
\def\radius{0.09}

\usepackage{matlab-prettifier}
\lstdefinestyle{matlab-custom}{
	language=Matlab,
	basicstyle=\footnotesize\ttfamily,
	keywordstyle=\bfseries\color{green!40!black},
	commentstyle=\itshape\color{purple!40!black},
	identifierstyle=\color{blue},
	stringstyle=\color{orange}
}

\makeatletter
\newcommand{\multiline}[1]{%
  \begin{tabularx}{\dimexpr\linewidth-\ALG@thistlm}[t]{@{}X@{}}
    #1
  \end{tabularx}
}
\makeatother

\usetikzlibrary{decorations}
\makeatletter
\def\pgf@dec@dashon{5pt}
\def\pgf@dec@dashoff{5pt}
\pgfkeys{/pgf/decoration/.cd,
  dash on/.store in=\pgf@dec@dashon,
  dash off/.store in=\pgf@dec@dashoff
}
\pgfdeclaredecoration{aligned dash}{start}{
\state{start}[width=0pt, next state=pre-corner,persistent precomputation={
\pgfextract@process\pgffirstpoint{\pgfpointdecoratedinputsegmentfirst}%
\pgfextract@process\pgfsecondpoint{\pgfpointdecoratedinputsegmentlast}%
\pgfmathsetlengthmacro\pgf@dec@dashon{\pgf@dec@dashon}%
\pgfmathsetlengthmacro\pgf@dec@dashoff{\pgf@dec@dashoff}%
\pgfmathsetlengthmacro\pgf@dec@halfdash{\pgf@dec@dashon/2}%
}]{}
\state{pre-corner}[width=\pgfdecoratedinputsegmentlength, next state=post-corner, persistent precomputation={
  \pgfmathparse{int(ceil((\pgfdecoratedinputsegmentlength-\pgf@dec@dashon-\pgf@dec@dashoff)/(\pgf@dec@dashon+\pgf@dec@dashoff)))}%
  \let\pgf@n=\pgfmathresult
  \pgfmathsetlengthmacro\pgf@b%
    {\pgfdecoratedinputsegmentlength/(\pgf@n+1)-\pgf@dec@dashon}%
  \ifdim\pgf@b<\pgf@dec@dashoff\relax%
    \pgfmathparse{int(\pgf@n-1)}\let\pgf@n=\pgfmathresult%
    \pgfmathsetlengthmacro\pgf@b%
      {\pgfdecoratedinputsegmentlength/(\pgf@n+1)-\pgf@dec@dashon}%
  \fi%
  \pgfmathsetlengthmacro\pgf@b{\pgf@b+\pgf@dec@dashon}%
}]{%
  \pgfmathloop
  \ifnum\pgfmathcounter>\pgf@n%
  \else%
    \pgfpathmoveto{\pgfpoint{\pgf@b*\pgfmathcounter-\pgf@dec@halfdash}{0pt}}%
    \pgfpathlineto{\pgfpoint{\pgf@b*\pgfmathcounter+\pgf@dec@halfdash}{0pt}}%
  \repeatpgfmathloop%
  \pgfpathmoveto%
    {\pgfpoint{\pgfdecoratedinputsegmentlength-\pgf@dec@halfdash}{0pt}}%
  \pgfpathlineto%
    {\pgfpointdecoratedinputsegmentlast}
}
\state{post-corner}[width=0pt, next state=pre-corner]{
   \pgfpathlineto{\pgfpoint{\pgf@dec@halfdash}{0pt}}%
}
\state{final}{
  \pgftransformreset%
  \pgfpathlineto{\pgfpointlineatdistance{\pgf@dec@halfdash}{\pgffirstpoint}{\pgfsecondpoint}}%
}
}
\tikzset{aligned dash/.style={
  decoration={aligned dash, #1}, decorate
}}

\def\ultraSEM{\texttt{\upshape ultraSEM}\xspace}
\newcommand{\dtn}{\Sigma}

\renewcommand{\vec}[1]{\bm{#1}}

\newcommand{\mytt}[1]{{\mlttfamily #1}}

\graphicspath{{figures/}}
\makeatletter
\def\input@path{{figures/}}
\makeatother

\begin{document}

\title{The ultraspherical spectral element method\thanks{Submitted to the editors December 15, 2020.
\funding{The first author was supported by the Department of Defense (DoD) through the National Defense Science \& Engineering Graduate Fellowship (NDSEG) Program. The second author was supported by the National Research Foundation (NRF) of South Africa (grant 109210). The third author was supported by NSF grant 1818757.}}}
\author{Daniel Fortunato\thanks{School of Engineering and Applied Sciences, Harvard University, Cambridge, MA 02138 (\email{dfortunato@g.harvard.edu}).} \and Nicholas Hale\thanks{Department of Mathematical Sciences, Stellenbosch University, Stellenbosch, 7602, South Africa (\email{nickhale@sun.ac.za}).} \and Alex Townsend\thanks{Department of Mathematics, Cornell University, Ithaca, NY 14853 (\email{townsend@cornell.edu}).}}
\headers{The ultraspherical spectral element method}{Daniel Fortunato, Nicholas Hale, and Alex Townsend}

\maketitle

\begin{abstract}
We introduce a novel spectral element method based on the ultraspherical spectral method and the hierarchical Poincar\'{e}--Steklov scheme for solving second-order linear partial differential equations on polygonal domains with unstructured quadrilateral or triangular meshes. Properties of the ultraspherical spectral method lead to almost banded linear systems, allowing the element method to be competitive in the high-polynomial regime ($p > 5$). The hierarchical Poincar\'{e}--Steklov scheme enables precomputed solution operators to be reused, allowing for fast elliptic solves in implicit and semi-implicit time-steppers. The resulting spectral element method achieves an overall computational complexity of $\mathcal{O}(p^4/h^3)$ for mesh size $h$ and polynomial order $p$, enabling $hp$-adaptivity to be efficiently performed. We develop an open-source software system, \ultraSEM, for flexible, user-friendly spectral element computations in MATLAB.
\end{abstract}



\section{Introduction}\label{sec:introduction}

Traditional approaches for solving partial differential equations (PDEs) on meshed geometries include finite element methods (FEMs)~\cite{Hughes_12_01}, discontinuous Galerkin (DG) methods~\cite{Cockburn_00_01}, and spectral element methods (SEMs)~\cite{Patera_84_01}. Each approach typically represents the solution of the PDE as a piecewise polynomial, with continuity or jump conditions weakly or strongly imposed between elements. Convergence is achieved by either refining the mesh ($h$-refinement) or increasing the polynomial degree on the elements ($p$-refinement). In theory, super-algebraic convergence can be observed---even for solutions with singularities---by optimally selecting a refinement strategy ($hp$-adaptivity)~\cite{Babuska_86_01}. However, $hp$-adaptivity theory can require high polynomial degrees, which are rarely used in practice as traditional methods can have prohibitive computational costs and numerical stability issues in this regime.

In particular, constructing efficient solvers for traditional high-order nodal element methods can be challenging. Direct solvers can become computationally intractable even for relatively small polynomial degrees as nodal discretizations result in dense linear algebra; in $d$ dimensions, the computational complexity for a direct solver na\"{i}vely scales as $\mathcal{O}(p^{3d})$. Iterative solvers may require an increasing number of iterations as $p$ increases because of the difficulties in designing robust preconditioners in the high $p$ regime~\cite{Orszag_80_01}. Because of these challenges, traditional element methods are typically restricted to low polynomial degrees, and $h$-refinement is generically preferred over $p$-refinement irrespective of local error estimators~\cite{Vos_10_01}. In practice, the physical considerations of the PDE---informed by $hp$-adaptivity theory---can take a back seat to the practical considerations of the numerical method.

Much work has gone toward reducing the computational costs associated with high-order element methods. For discretizations that possess tensor-product structure (e.g., standard nodal bases on quadrilateral elements or certain bases on triangular elements~\cite{Sherwin_95_01}), sum factorization~\cite{Orszag_80_01} reduces the cost of operator assembly from $\mathcal{O}(p^{3d})$ to $\mathcal{O}(p^{2d+1})$, and matrix-free evaluation reduces the cost of matrix-vector multiplication from $\mathcal{O}(p^{2d})$ to $\mathcal{O}(p^{d+1})$~\cite[Tab.~1]{MFEM}. Solvers for the resulting linear systems are often based on iterative methods coupled with sufficient preconditioning. Low-order FEM discretizations on a mesh constructed from the high-order SEM nodes can be shown to be spectrally equivalent to the SEM discretizations~\cite{Canuto_07_01}, and matrix-free preconditioners based on this equivalence can perform well when coupled with a multigrid method using specialized smoothers~\cite{Pazner_19_01}. Multigrid methods applied to high-order DG discretizations can perform well if the discrete operators are coarsened according to the flux formulation of the PDE~\cite{Fortunato_19_02}. Spectral element multigrid methods have proven effective when applied to nodal discretizations of Poisson's equation in one dimension, though multigrid convergence factors can weakly depend on both $h$ and $p$~\cite{Ronquist_87_01,Maday_88_01}. Modal discretizations for $p$-FEM based on integrated Jacobi polynomials can yield sparse stiffness matrices that contain an optimal number of nonzeros, but developing optimal solvers for such discretizations remains a challenge~\cite{Beuchler_06_01}. Many open-source software libraries exist for high-order element computation, including MFEM~\cite{MFEM}, Firedrake~\cite{Firedrake}, Nektar++~\cite{Nektar,Nektar2}, and Nek5000~\cite{Nek5000}.

Though solvers for element methods are commonly based on preconditioned iterative methods, fast direct solvers for high-order methods have become an active area of research in recent years. The hierarchical Poincar\'{e}--Steklov (HPS) scheme~\cite{Martinsson_09_01, Martinsson_13_01, Gillman_14_01, Gillman_14_02, Babb_18_01} is a multidomain spectral collocation method based on a recursive domain decomposition approach, which ``glues'' together solutions at interfaces between elements using Poincar\'{e}--Steklov operators (such as Dirichlet-to-Neumann operators). The accompanying direct solver is analagous to classical nested dissection. The formulation hierarchically merges Dirichlet-to-Neumann operators and results in an in-memory solution operator, which can be reapplied fast to multiple righthand sides on static meshes. The ability to reuse computed solution operators allows for efficient implicit time-stepping for parabolic problems~\cite{Babb_18_02}. The method has been extended to handle mesh adaptivity~\cite{Geldermans_19_01}, three-dimensional problems~\cite{Hao_16_01}, and boundary integral equations~\cite{Gillman_15_01}. The HPS scheme based on spectral collocation has an overall complexity of $\mathcal{O}(Np^4 + N^{3/2})$, where $N \approx (p/h)^2$ is the total number of degrees of freedom, $p$ is the polynomial degree on each element, and $h$ is the minimum mesh element size.

In this paper, we take advantage of recent advances in sparse spectral methods to propose an SEM in two dimensions with a computational complexity of
\vspace{0.2em}
\[
\underbrace{\frac{p^4}{h^2} \,+\, \frac{p^3}{h^3}}_{\text{build stage}} \;\;+\;\; \underbrace{\frac{p^3}{h^2} \,+\, \frac{p^2}{h^2}\log \frac{1}{h^2}}_{\text{solve stage}} \;\approx\; \frac{p^4}{h^2} + \frac{p^3}{h^3} \;\approx\; Np^2 + N^{3/2}.\vspace{0.2em}
\]
Specifically, we propose a variant of the HPS scheme that employs the ultraspherical spectral method~\cite{Olver_13_01,Townsend_15_01} instead of spectral collocation for element-wise discretization. The method retains sparsity in the high-$p$ regime by carefully selecting bases to be specific families of orthogonal polynomials and employing sparse recurrence relations between them. The discretization is not nodal, but modal; that is, the unknowns are not values on a grid, but coefficients in a polynomial expansion.

In this work, we are interested in solving linear PDEs on two-dimensional meshed geometries with Dirichlet boundary conditions,\footnote{Robin boundary conditions can be converted to equivalent Dirichlet boundary conditions using the Dirichlet-to-Neumann operators constructed by the HPS scheme, and so we focus on Dirichlet boundary conditions throughout the paper.
} i.e., 
\begin{equation}\label{eq:pde}
\begin{aligned}
\mathcal{L}u(x,y) &= f(x,y) &&\text{in } \Omega, \\
u(x,y) &= g(x,y) &&\text{on } \partial\Omega.
\end{aligned}
\end{equation}
Here, $\Omega$ is a domain in $\mathbb{R}^2$, $f$ and $g$ are given functions defined on $\Omega$ and its boundary, and $\mathcal{L}$ is a variable-coefficient, second-order, elliptic partial differential operator (PDO) of the form
\begin{equation}\label{eq:pdo}
\mathcal{L}u = \nabla \cdot \left(A(x,y) \nabla u\right) + \nabla \cdot \left(b(x,y) u\right) + c(x,y) u,
\end{equation}
with $A(x,y) \in \mathbb{C}^{2 \times 2}$, $b(x,y) \in \mathbb{C}^2$, and $c(x,y) \in \mathbb{C}$.

The paper is structured as follows. In \cref{sec:background}, we review the ultraspherical spectral method, a sparse and spectrally-accurate method for solving linear ODEs and PDEs on rectangular domains, and discuss its application to quadrilateral and triangular domains. In \cref{sec:ultraSEM}, we extend this spectral method to the non-overlapping domain decomposition setting, highlighting the differences from traditional collocation-based patching approaches. We describe how the hierarchical merging of Poincar\'{e}--Steklov operators efficiently performs domain decomposition on meshes with many elements. In \cref{sec:software}, we present an implementation of the ultraspherical SEM in the software package \ultraSEM, and briefly describe its syntax and design. In \cref{sec:results}, we present numerical results and applications of the method.

\section{Background material}\label{sec:background}

\subsection{The ultraspherical spectral method}\label{sec:ultraS}

First, we review the fundamental ideas in the ultraspherical spectral method~\cite{Olver_13_01}, which in one dimension solves linear ordinary differential equations (ODEs) with variable coefficients of the form
\begin{equation}\label{eq:ode}
\sum_{\lambda=0}^M a_\lambda(x)\frac{d^\lambda u}{dx^\lambda} = f(x), \qquad x\in[-1,1],
\end{equation}
along with general linear boundary conditions $\mathcal{B}u=\vec{\mathrm{g}} \in \mathbb{C}^M$ to ensure that there is a unique solution. For an integer $p$, the method seeks to approximate the first $p+1$ Chebyshev expansion coefficients $\{u_j\}_{j=0}^p$ of the solution $u$, where
\[
u(x) = \sum_{j=0}^\infty u_j T_j(x), \qquad x \in [-1,1],
\]
and $T_j(x) = \cos(j\cos^{-1}x)$ is the degree-$j$ Chebyshev polynomial of the first kind.

Classical spectral methods represent differentiation as a dense operator~\cite{Trefethen_00_01, Boyd_01_01}, but the ultraspherical spectral method employs the ``sparse'' recurrence relations
\begin{equation}
\frac{d^\lambda T_j}{dx^\lambda} = \begin{cases}
2^{\lambda-1} j (\lambda-1)!\,C^{(\lambda)}_{j-\lambda}, & j \geq \lambda, \\
0, & 0 \leq j \leq \lambda-1,
\end{cases}
\end{equation}
where $C^{(\lambda)}_j$ is the degree-$j$ ultraspherical polynomial of parameter $\lambda > 0$~\cite[Sec. 18.3]{NISTHandbook}. This results in a sparse representation of differentiation operators. In particular, the differentiation operator for the $\lambda$th derivative is given by
\[
\mathcal{D}_\lambda = 2^{\lambda-1} (\lambda-1)! \begin{pmatrix}
\overbrace{0\;\;\cdots\;\;0}^{\lambda \text{ times}} & \lambda \\
&& \lambda+1 \\
&&& \lambda+2 \\
&&&& \ddots
\end{pmatrix}, \qquad \lambda \in \mathbb{N}\setminus\{0\}.
\]
For $\lambda \geq 1$, the matrix $\mathcal{D}_\lambda$ maps a vector of Chebyshev coefficients to a vector of $C^{(\lambda)}$ coefficients of the $\lambda$th derivative. For convenience, we use $\mathcal{D}_0$ to denote the identity operator.

Since $\mathcal{D}_\lambda$ returns a vector of ultraspherical coefficients for $\lambda\geq1$, operators to convert between the Chebyshev and ultraspherical bases are required. Let $\mathcal{S}_0$ be the operator that converts a vector of Chebyshev coefficients to a vector of $C^{(1)}$ coefficients, and let $\mathcal{S}_\lambda$, for $\lambda \geq 1$, be the operator that converts a vector of $C^{(\lambda)}$ coefficients to a vector of $C^{(\lambda+1)}$ coefficients. Using the recurrence relations~\cite[(18.9.7) \& (18.9.9)]{NISTHandbook}
\[
T_j = \begin{cases}
\tfrac{1}{2}\left(C_j^{(1)}-C_{j-2}^{(1)}\right), & j\geq 2, \\[0.4em]
\tfrac{1}{2}C_1^{(1)}, & j = 1, \\[0.3em]
C_0^{(1)}, & j = 0,
\end{cases}
\qquad
C_j^{(\lambda)} = \begin{cases}
\tfrac{\lambda}{\lambda+j} \left(C_{j}^{(\lambda+1)}-C_{j-2}^{(\lambda+1)}\right), & j\geq 2, \\[0.4em]
\tfrac{\lambda}{\lambda+1} C_1^{(\lambda+1)}, & j = 1, \\[0.3em]
C_0^{(\lambda+1)}, & j = 0,
\end{cases}
\]
it can be shown that the conversion operators $\mathcal{S}_0$ and $\mathcal{S}_\lambda$ are sparse and given by~\cite{Olver_13_01}
\[
\mathcal{S}_0 = \begin{pmatrix}
1 & 0 & -\frac12 \\[0.5em]
& \frac12 & 0 & -\frac12 \\
&& \frac12 & 0 & \ddots \\
&&& \frac12 & \ddots \\
&&&& \ddots
\end{pmatrix}, \quad
\mathcal{S}_\lambda = \begin{pmatrix}
1 & 0 & -\frac{\lambda}{\lambda+2} \\[0.5em]
& \frac{\lambda}{\lambda+1} & 0 & -\frac{\lambda}{\lambda+3} \\
&& \frac{\lambda}{\lambda+2} & 0 & \ddots \\
&&& \frac{\lambda}{\lambda+2} & \ddots \\
&&&& \ddots
\end{pmatrix},
\quad \lambda \geq 1.
\]

To represent multiplication by the variable coefficients $a_\lambda(x)$ in \cref{eq:ode}, multiplication operators $\mathcal{M}_\lambda[a_\lambda]$ for $C^{(\lambda)}$ coefficients\footnote{The multiplication operator for $\lambda = 0$, $\mathcal{M}_0[a_0]$, acts on a vector of Chebyshev coefficients.} can be explicitly constructed. If $a_\lambda(x)$ is approximated by a degree-$m_\lambda$ polynomial, then the operator $\mathcal{M}_\lambda[a_\lambda]$ is $m_\lambda$-banded~\cite{Olver_13_01}.

Discretizing \cref{eq:ode} using these operators to represent differentiation, conversion between bases, and multiplication by variable coefficients results in a banded $(p+1) \times (p+1)$ linear system given by
\begin{equation}\label{eq:ode_ultra}
\left( \mathcal{M}_M[a_M] \mathcal{D}_M + \sum_{\lambda=0}^{M-1} \mathcal{S}_{M-1} \cdots \mathcal{S}_\lambda \mathcal{M}_\lambda[a_\lambda] \mathcal{D}_\lambda \right) \vec{\mathrm{u}} = \mathcal{S}_{M-1} \cdots \mathcal{S}_0 \, \vec{\mathrm{f}},
\end{equation}
where $\vec{\mathrm{u}}$ and $\vec{\mathrm{f}}$ are vectors of Chebyshev coefficients of $u$ and $f$, respectively. Note that since the order-$M$ differential operator in \cref{eq:ode_ultra} maps the vector of Chebyshev coefficients $\vec{\mathrm{u}}$ to $C^{(M)}$ coefficients, the vector of Chebyshev coefficients $\vec{\mathrm{f}}$ must also be converted to $C^{(M)}$ coefficients. The bandwidth of the linear system in \cref{eq:ode_ultra} scales as $\mathcal{O}(\max_\lambda m_\lambda)$, independent of the polynomial order $p$. If the variable coefficients $a_\lambda(x)$ can be approximated by polynomials such that $m_\lambda \ll p$, then \cref{eq:ode_ultra} is a sparse linear system.

\begin{figure}
  \centering
  \includegraphics[width=3.97cm,trim=0 0.3cm 0 0]{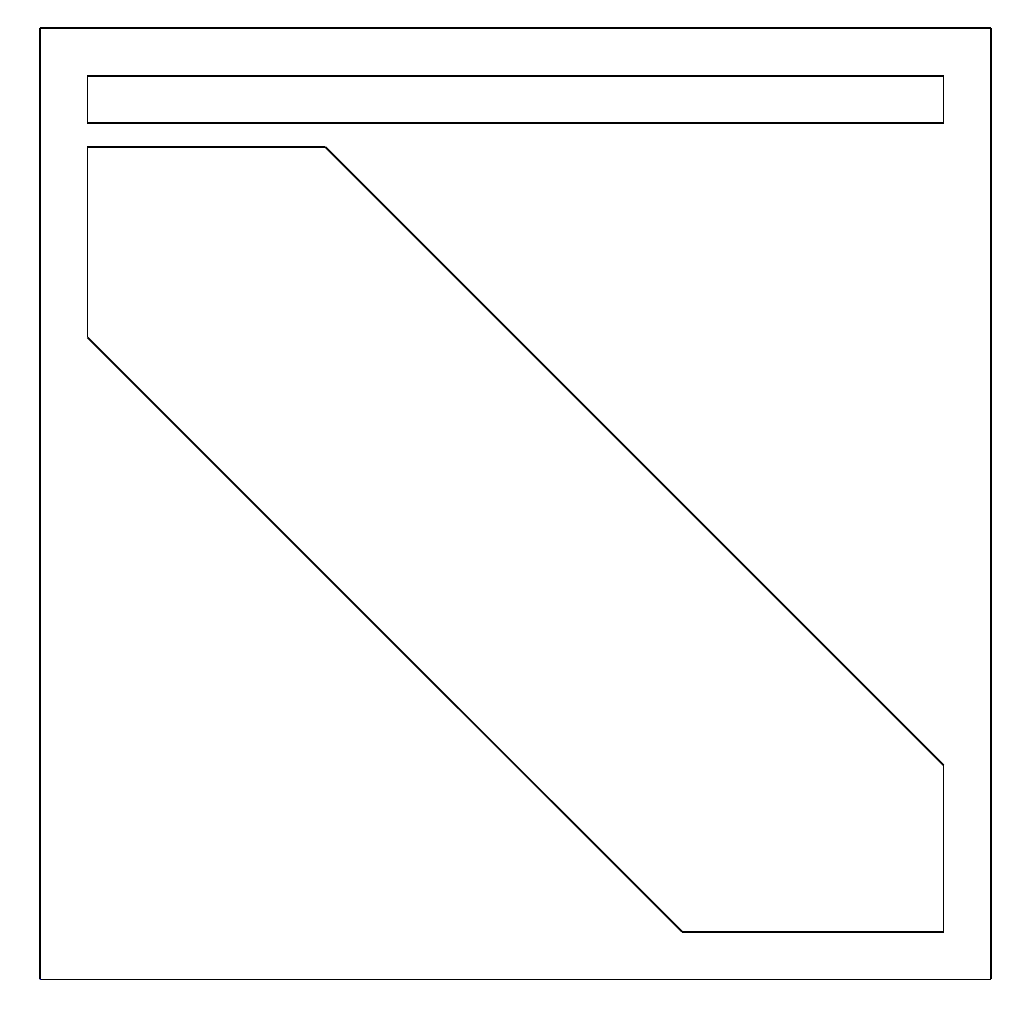}%
  \hspace{0.16cm}
  \includegraphics[width=3.7cm]{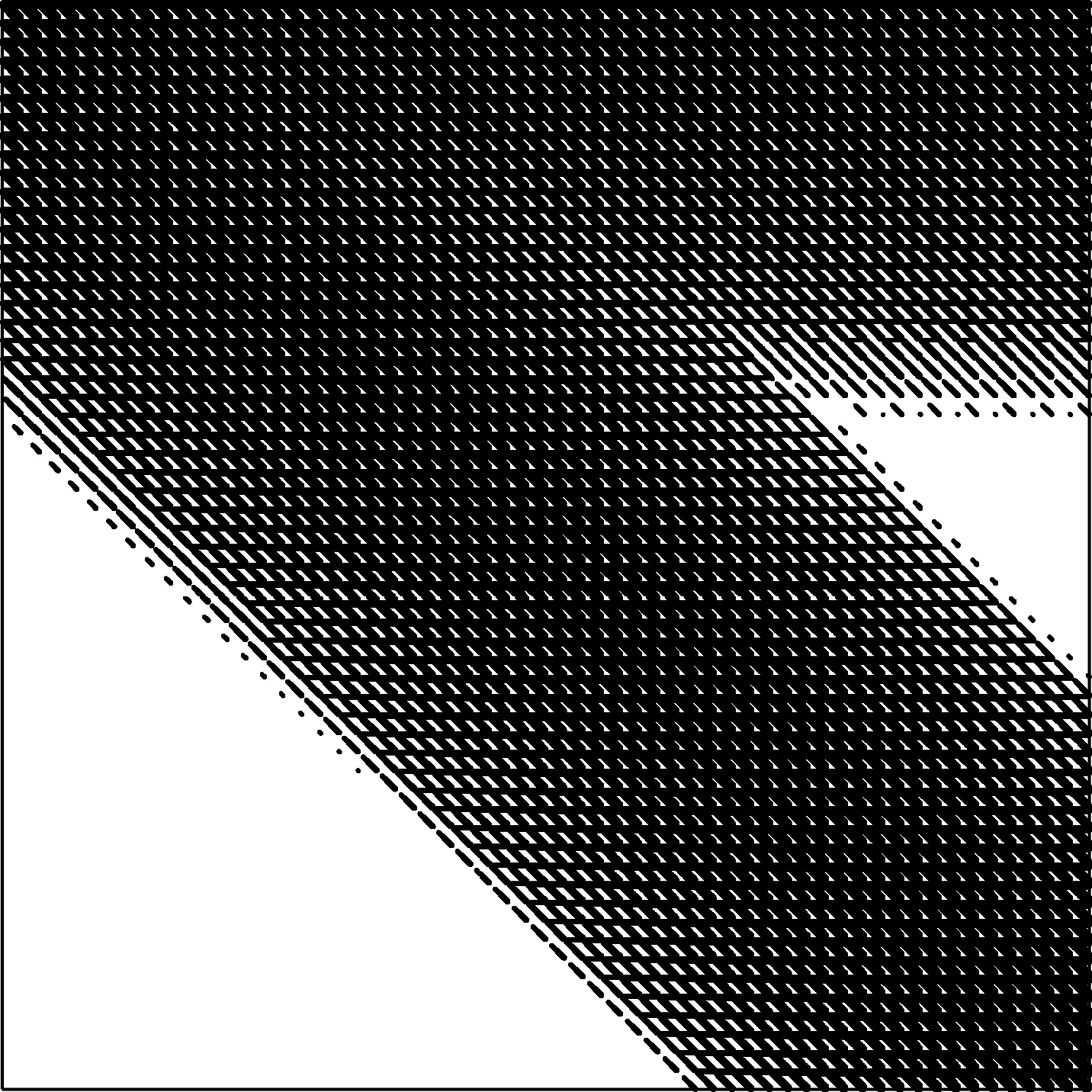}%
  \hspace{0.25cm}
  \includegraphics[width=3.7cm]{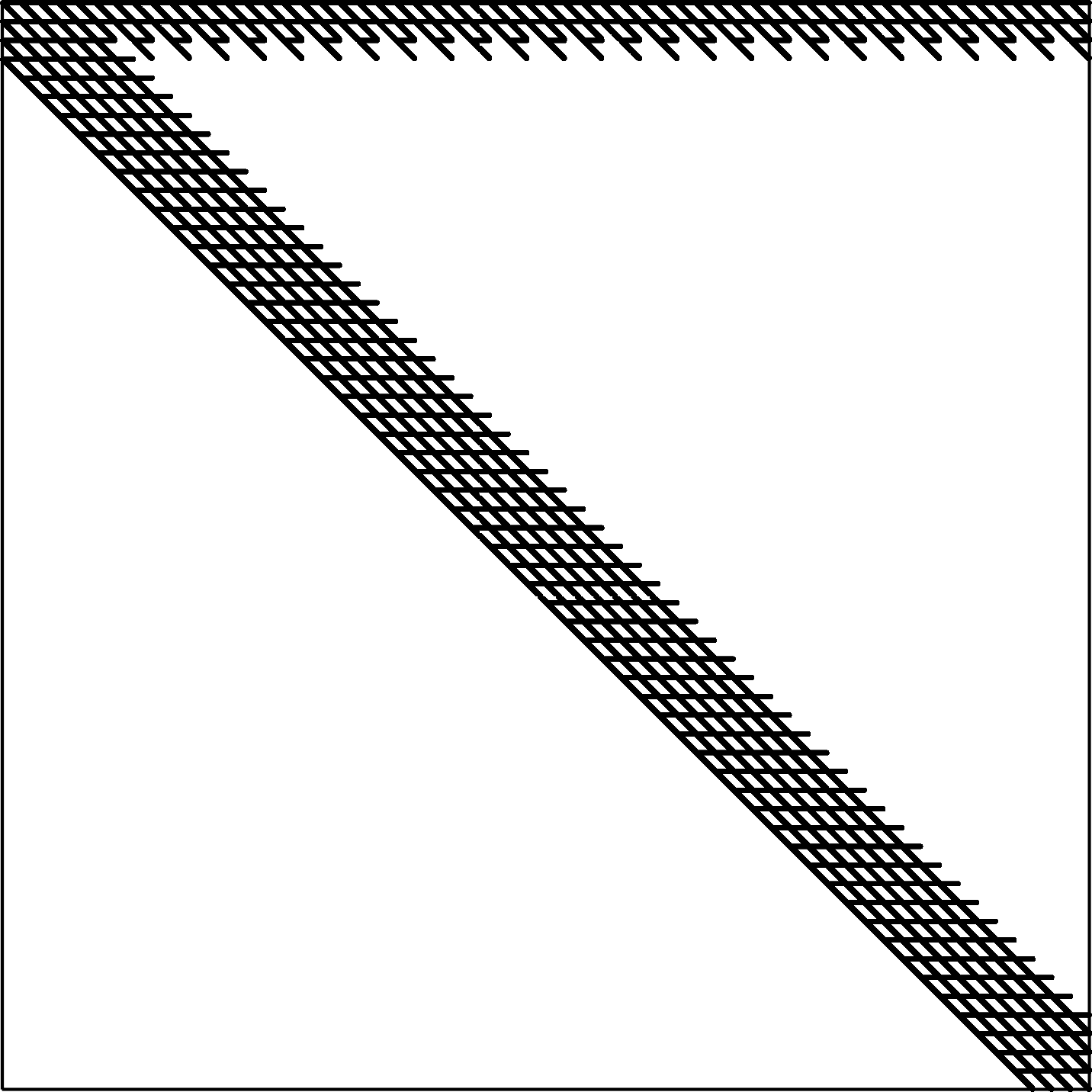}%
\caption{(Left) Typical structure of the almost banded matrices constructed by the ultraspherical spectral method, i.e., banded matrices except for a small number of dense rows. In one dimension, ODEs are discretized as almost banded $(p+1) \times (p+1)$ linear systems with bandwidth $m$ independent of $p$ and $\mathcal{O}(1)$ dense rows; such systems can be solved in $\mathcal{O}(p)$ operations. In two dimensions on $[-1,1]^2$, PDEs are discretized as almost block-banded $(p+1)^2 \times (p+1)^2$ linear systems, with a bandwidth of $\mathcal{O}(m p)$ and $\mathcal{O}(m p)$ dense rows; such systems can be solved in $\mathcal{O}(p^4)$ operations.
(Center) In two dimensions on a quadrilateral, PDEs are transformed to $[-1,1]^2$ and then discretized. When Jacobian factors are kept as rational functions (see \cref{sec:quad_tri}), the discrete differential operator has a large bandwidth. (Right) By scaling the transformed PDE by a power of the determinant of the Jacobian, the discrete differential operator remains sparse. This improvement is most notable when the quadrilaterals involved have small angles.}
\label{fig:bandedness}
\end{figure}

To impose the boundary constraints given by $\mathcal{B}$, we must encode $\mathcal{B}$ in terms of its action on a vector of Chebyshev coefficients. For Dirichlet boundary conditions on $[-1,1]$, such action is given by
\begingroup
\renewcommand*{\arraystretch}{1.2}
\begin{equation}\label{eq:dir_bc}
\mathcal{B} = \begin{pmatrix}
T_0(-1) & T_1(-1) & \cdots & T_p(-1) \\
T_0(1) & T_1(1) & \cdots & T_p(1)
\end{pmatrix} = \begin{pmatrix}
1 & -1 & \cdots & (-1)^p \\
1 & 1 & \cdots & 1
\end{pmatrix},
\end{equation}
\endgroup
because $\mathcal{B}\vec{\mathrm{u}} \approx (u(-1), u(1))^T$. Neumann, Robin, and more general boundary constraints can be similarly encoded. To impose the $M$ boundary conditions $\mathcal{B}\vec{\mathrm{u}}=\vec{\mathrm{g}}$ on the linear system \cref{eq:ode_ultra}, the ultraspherical spectral method uses boundary bordering~\cite{Boyd_01_01}, wherein the last $M$ rows of the linear system are replaced by dense rows that impose constraints on the Chebyshev coefficients of the solution (e.g., \cref{eq:dir_bc} for Dirichlet boundary conditions). The resulting $(p+1) \times (p+1)$ linear system has a distinctive almost banded\footnote{A matrix is \emph{almost banded} if it is banded except for a small number of columns or rows.} structure with bandwidth $m = \max_\lambda m_\lambda$ and can be solved in $\mathcal{O}(m^2 p)$ operations using the adaptive QR algorithm~\cite{Olver_13_01} or the Woodbury formula. \cref{fig:bandedness} (left) shows the almost banded structure typical of the linear systems in the ultraspherical spectral method.

The ultraspherical spectral method can be extended to solve PDEs in two dimensions on rectangular domains~\cite{Townsend_15_01}. For the PDE given in \cref{eq:pde} and for a polynomial order $p$, the method computes modes $X \in \mathbb{C}^{(p+1) \times (p+1)}$ of the solution $u(x,y)$ in a bivariate tensor-product Chebyshev basis, such that
\[
u(x,y) = \sum_{i=0}^\infty \sum_{j=0}^\infty X_{ij} T_i(y) T_j(x), \qquad (x,y) \in [-1,1]^2.
\]
Discretization of the PDE is based on separable models of linear partial differential operators. For example, the elliptic PDO $\mathcal{L}$ given by \cref{eq:pdo} can be decomposed into a sum of tensor products of one-dimensional differential operators
\begin{equation}
\mathcal{L} = \sum_{j=1}^K \left(\mathcal{L}_j^y\otimes \mathcal{L}_j^x\right),
\label{eq:separableModel}
\end{equation}
where $\mathcal{L}^y_1,\ldots, \mathcal{L}^y_K$ are operators associated with ODEs in $y$, $\mathcal{L}_1^x,\ldots, \mathcal{L}_K^x$ are operators associated with ODEs in $x$. In~\cref{eq:separableModel}, the tensor product operator `$\otimes$' is defined such that if $u(x,y) = v(y) w(x)$, then
\[
\left(\mathcal{L}^y \otimes \mathcal{L}^x\right) u(x,y) = \left(\mathcal{L}^y v(y)\right) \left(\mathcal{L}^x w(x)\right)
\]
for some operators $\mathcal{L}^y$ and $\mathcal{L}^x$.
Such separable representations of PDOs can be automatically computed~\cite{Townsend_15_01}. The univariate differential operators $\mathcal{L}^y_1,\ldots, \mathcal{L}^y_K, \mathcal{L}_1^x,\ldots, \mathcal{L}_K^x$ can each be discretized using the ultraspherical spectral method in one dimension, and boundary conditions in $x$ and $y$ can be imposed on the rows and columns of $X$, thus giving us a scheme for discretizing the PDE. The resulting linear system of size $(p+1)^2 \times (p+1)^2$ is almost block-banded with a bandwidth of $\mathcal{O}(m p)$ and $\mathcal{O}(m p)$ dense rows, where $m = \max\{m_1^y, \ldots, m_K^y, m_1^x, \ldots, m_K^x\}$ and $m_1^y, \ldots, m_K^y, m_1^x, \ldots, m_K^x$ are the bandwidths of the discretized operators $\mathcal{L}^y_1,\ldots, \mathcal{L}^y_K, \mathcal{L}_1^x,\ldots, \mathcal{L}_K^x$. This can be solved in $\mathcal{O}(p^4)$ operations. In special cases, e.g., where $K=1$ or $K=2$, further structure can be exploited to arrive at faster solvers~\cite{Townsend_15_01,Fortunato_19_01}.

\subsection{Spectral methods on quadrilaterals and triangles}\label{sec:quad_tri}
Global spectral methods defined on rectangles can be used on other polygons through coordinate transformation. Let $\mathcal{Q}_\text{ref} = [-1,1]^2$ be the reference square with vertices given by $(r_0,s_0) = (-1,-1)$, $(r_1,s_1) = (1,-1)$, $(r_2,s_2) = (1,1)$, $(r_3,s_3) = (-1,1)$. Denote by $(r,s)$ the coordinates in reference space and by $(x,y)$ the coordinates in real space, and suppose we have a mapping from reference space to real space, $(r,s) \mapsto (x,y)$. To apply a global spectral method on $\mathcal{Q}_\text{ref}$ to a PDE defined in real space, the differential operator $\mathcal{L}$ and righthand side $f(x,y)$ are transformed into reference space. The coordinate transformation alters the differential operator via the chain rule. For a function $u(r,s)$ defined on $\mathcal{Q}_\text{ref}$, first- and second-order derivatives in $x$ and $y$ are given by
\[
\begin{aligned}
u_x &= r_x u_r + s_x u_s, \\
u_y &= r_y u_r + s_y u_s, \\
u_{xx} &= (r_x)^2 u_{rr} + 2 r_x s_x u_{rs} + (s_x)^2 u_{ss} + r_{xx} u_r + s_{xx} u_s, \\
u_{xy} &= r_x r_y u_{rr} + (r_x s_y + r_y s_x) u_{rs} + s_x s_y u_{ss} + r_{xy} u_r + s_{xy} u_s, \\
u_{yy} &= (r_y)^2 u_{rr} + 2 r_y s_y u_{rs} + (s_y)^2 u_{ss} + r_{yy} u_r + s_{yy} u_s, \\
\end{aligned}
\]
where the Jacobian factors $r_x, r_{xx}, \ldots$ depend on the coordinate mapping. In this paper, we are interested in mappings from $\mathcal{Q}_\text{ref}$ to quadrilaterals or triangles.

For a quadrilateral domain $\mathcal{Q}$ with vertices $(x_0,y_0), \ldots, (x_3,y_3)$, a bilinear mapping from $(r,s) \in \mathcal{Q}_\text{ref}$ to $(x,y) \in \mathcal{Q}$ is given by
\[
\begin{bmatrix} r \\ s \end{bmatrix}
\mapsto
\begin{bmatrix} a_0^x + a_1^x r + a_2^x s + a_3^x r s \\[0.2em] a_0^y + a_1^y r + a_2^y s + a_3^y r s \end{bmatrix}
=
\begin{bmatrix} x \\ y \end{bmatrix},
\]
where the coefficients $a_0^x, \ldots, a_3^x$ and $a_0^y, \ldots, a_3^y$ satisfy the linear system

\[
\begin{bmatrix}
1 & r_0 & s_0 & r_0 s_0 \\
1 & r_1 & s_1 & r_1 s_1 \\
1 & r_2 & s_2 & r_2 s_2 \\
1 & r_3 & s_3 & r_3 s_3
\end{bmatrix}
\begin{bmatrix}
a_0^x & a_0^y \\
a_1^x & a_1^y \\
a_2^x & a_2^y \\
a_3^x & a_3^y
\end{bmatrix}
=
\begin{bmatrix}
x_0 & y_0 \\
x_1 & y_1 \\
x_2 & y_2 \\
x_3 & y_3
\end{bmatrix}.
\]
While the mapping from $(r,s)$ to $(x,y)$ is bilinear, the mapping from $(x,y)$ to $(r,s)$ is more complicated and in particular is not polynomial, and so we would like to avoid directly computing the inverse maps $r(x,y)$ and $s(x,y)$. Therefore, to compute the first-order Jacobian factors $r_x, s_x, r_y$, and $s_y$, we apply the inverse function theorem to the Jacobian matrix $J_{rs} = \partial (r,s) / \partial (x,y)$, which states that $J_{rs} = \left(J_{xy}\right)^{-1}$ with $J_{xy} = \partial (x,y) / \partial (r,s)$. Writing out the Jacobians explicitly, we obtain the following formulae for the first-order factors $r_x, s_x, r_y$, and $s_y$:
\[
\begin{bmatrix}
r_x & r_y \\
s_x & s_y
\end{bmatrix}
=
\begin{bmatrix}
x_r & x_s \\
y_r & y_s
\end{bmatrix}^{-1}
=
\frac{1}{\det(J_{xy})}
\begin{bmatrix}
\phantom{-}y_s & -x_s \\
-y_r & \phantom{-}x_r
\end{bmatrix},
\]
where $\det(J_{xy}) = x_r y_s - x_s y_r$. Applying the chain rule to these definitions yields formulae for the second-order factors $r_{xx}$, $r_{xy}$, $r_{yy}$, $s_{xx}$, $s_{xy}$, and $s_{yy}$.

However, note that the Jacobian factors are rational functions, due to factors of $\det(J_{xy})$, $\det(J_{xy})^2$, and $\det(J_{xy})^3$ in the denominators of the first- and second-order terms. Thus, the coordinate transformation from $\mathcal{Q}$ to $\mathcal{Q}_\text{ref}$ introduces rational variable coefficients into the differential operator, and the discretization of the transformed operator by the ultraspherical spectral method results in a linear system with large bandwidth (see \cref{fig:bandedness} (center)). To recover sparsity, we scale the transformed differential operator $\mathcal{L}_{rs}$ and righthand side $f(r,s)$ by the factor $\det(J_{xy})^3$~\cite{Yeiser_18_01}, and discretize the scaled PDE
\[
\underbrace{\left( \det(J_{xy})^3 \mathcal{L}_{rs} \right)}_{\widehat{\mathcal{L}}_{rs}} u(r,s) = \underbrace{ \det(J_{xy})^3 f(r,s)}_{\widehat{f}}.
\vspace{-0.3em}
\]
As all Jacobian factors can be written with denominator $\det(J_{xy})^3$, this scaling turns the rational variable coefficients induced by the transformation into polynomial variable coefficients of degree $\leq 3$ (see \cref{fig:bandedness} (right)). Thus, PDEs on $\mathcal{Q}$ with degree-$m$ variable coefficients are transformed into PDEs on $\mathcal{Q}_\text{ref}$ with degree-$(m+3)$ variable coefficients.

For a triangular domain $\mathcal{T}$, the Duffy transformation~\cite{Duffy_82_01,Sherwin_05_01} may be used to define a mapping from $\mathcal{Q}_\text{ref}$ to $\mathcal{T}$ by collapsing one side of $\mathcal{Q}_\text{ref}$ to a point. Let $\mathcal{T}_\text{ref}$ be the reference triangle with vertices $(x_0,y_0) = (0,0)$, $(x_1,y_1) = (1,0)$, and $(x_2,y_2) = (0,1)$. A mapping from $(r,s) \in \mathcal{Q}_\text{ref}$ to $(x,y) \in \mathcal{T}_\text{ref}$ can be defined by
\[
\begin{bmatrix} r \\ s \end{bmatrix}
\mapsto
\begin{bmatrix} \tfrac14(1+r)(1-s) \\[0.2em] \tfrac12(1-s) \end{bmatrix}
=
\begin{bmatrix} x \\ y \end{bmatrix},
\]
which maps the line segment between $(-1,1)$ and $(1,1)$ in $\mathcal{Q}_\text{ref}$ to the point $(0,1)$ in $\mathcal{T}_\text{ref}$. The inverse of this transformation, mapping from $(x,y) \in \mathcal{T}_\text{ref}$ to $(r,s) \in \mathcal{Q}_\text{ref}$, possesses a singularity at the point $(0,1)$, i.e.,
\[
\begin{bmatrix} x \\ y \end{bmatrix}
\mapsto
\begin{bmatrix} 2x/(1-y)-1 \\[0.2em] 2y-1 \end{bmatrix}
=
\begin{bmatrix} r \\ s \end{bmatrix}.
\]
If discretized directly, Jacobian factors based on this transformation introduce singular variable coefficients into the differential operator when the operator is transformed to $\mathcal{T}_\text{ref}$. However, the singularity induced by the Duffy transformation may be removed by scaling the PDE by powers of $1-y$. For a general triangular domain $\mathcal{T}$ with vertices $(x_0,y_0), \ldots, (x_2,y_2)$, the Duffy transformation may be composed with an affine transformation of the form
\[
\begin{bmatrix} x \\ y \end{bmatrix}
\mapsto
\begin{bmatrix}
x_0 + (x_1-x_0)\,x + (x_2-x_0)\,y \\[0.2em]
y_0 + (y_1-y_0)\,x + (y_2-y_0)\,y \end{bmatrix}
\]
to yield a mapping from $\mathcal{Q}_\text{ref}$ to $\mathcal{T}$.

We focus our attention on straight-sided quadrilateral elements in the remainder of this work. However, the algorithms presented below can be applied to triangular elements through simple modifications. The \ultraSEM software supports both triangular and quadrilateral elements.

\section{The ultraspherical spectral element method}\label{sec:ultraSEM}
We now describe how to adapt the ultraspherical spectral method into an SEM, focusing on key implementation aspects. Our method is based on the hierarchical Poincar\'{e}--Steklov scheme, an efficient non-overlapping domain decomposition approach~\cite{Martinsson_09_01, Martinsson_13_01, Gillman_14_01, Gillman_14_02, Babb_18_01}. We employ a variant of the HPS scheme to handle irregular, non-tensor-product meshes (see \cref{sec:complexity}). Broadly, our method is the following:

\vspace{1em}
\begin{enumerate}[itemsep=0.5em]

\item The method takes as input a second-order elliptic PDO $\mathcal{L}$, a righthand side $f$, Dirichlet data $g$, and a mesh with elements $\{\mathcal{E}_i\}_{i=1}^{n_\text{elem}}$.

\item On each element, two local operators are constructed: (i) a solution operator, which computes the local solution to the PDE on the element when given Dirichlet data, and (ii) a Dirichlet-to-Neumann operator, which computes the outward flux of the local solution when given Dirichlet data (see \cref{sec:ultra_hps:build}).

\item Local elemental operators are merged pairwise in a hierarchical fashion, yielding solution operators and Dirichlet-to-Neumann operators, which act on the interfaces between elements or groups of elements. Merging continues until a single global solution operator is computed for the entire mesh (see \cref{sec:ultra_hps:build}).

\item The given Dirichlet data $g$ is passed in at the top level. Solution operators are applied down the tree, providing the solution at unknown interfaces between elements (see \cref{sec:ultra_hps:solve}).

\item Once the solution is known at all the interfaces, local solution operators are applied on each element to determine the interior solution over the entire mesh.

\end{enumerate}
\vspace{1em}

The method naturally lends itself to parallelization. Specifically, steps 2 and 5 can be performed independently on each element as the computations involved are entirely decoupled. Moreover, step 2 is often the bottleneck when $p$ is large, and so significant speedups may be gained if parallelism is exploited (see \cref{sec:complexity}). The hierarchical steps 3 and 4 may also be parallelized, as the operations taking place on two branches in the hierarchy are decoupled until the two branches are merged. Thus, a careful, load-balanced strategy for parallelizing across branches in the hierarchy may lead to further speedups.

\subsection{Domain decomposition for modal discretizations}\label{sec:modal_dd}
\begin{figure}[htb]
  \centering
  \vspace{-0.3cm}
	\subfloat[Nodal discretization]{%
		\label{fig:vals_vs_coeffs:nodal}
\raisebox{0.31cm}{%
      \begin{tikzpicture}[scale=1]
      \begin{scope}[xshift=-0.0cm]
          \filldraw[fill=white, draw=black, very thick] (0,0) rectangle (2,2);
          \foreach \x in {0.133974596215561,0.5,1,1.5,1.866025403784439}
	        \foreach \y in {0.133974596215561,0.5,1,1.5,1.866025403784439}
              \fill[\intcolor] (\x,\y) circle (\radius);
          \foreach \x in {0,0.133974596215561,0.5,1,1.5,1.866025403784439,2}
          {
            \fill[\extcolor] (\x,0) circle (\radius);
            \fill[\extcolor] (\x,2) circle (\radius);
            \fill[\extcolor] (0,\x) circle (\radius);
            \fill[\extcolor] (2,\x) circle (\radius);
          }
      \end{scope}
      \begin{scope}[xshift=2cm]
          \filldraw[fill=white, draw=black, very thick] (0,0) rectangle (2,2);
          \foreach \x in {0.133974596215561,0.5,1,1.5,1.866025403784439}
	        \foreach \y in {0.133974596215561,0.5,1,1.5,1.866025403784439}
              \fill[\intcolor] (\x,\y) circle (\radius);
          \foreach \x in {0,0.133974596215561,0.5,1,1.5,1.866025403784439,2}
          {
            \fill[\extcolor] (\x,0) circle (\radius);
            \fill[\extcolor] (\x,2) circle (\radius);
            \fill[\extcolor] (0,\x) circle (\radius);
            \fill[\extcolor] (2,\x) circle (\radius);
          }
      \end{scope}
    \end{tikzpicture}}}
  \hspace{1cm}
	\subfloat[Modal discretization]{%
	\label{fig:vals_vs_coeffs:modal}
      \raisebox{0.15cm}{\begin{tikzpicture}[scale=1]
      \begin{scope}[xshift=-0.0cm]
          \filldraw[fill=white, draw=black, very thick, aligned dash={dash on=4pt, dash off=3pt}] (0,0) rectangle (2,2);
          \fill[fill=\extcolor] (-0.25,0) rectangle (-0.15,2);
          \fill[fill=\extcolor] (2.25,0) rectangle (2.15,2);
          \fill[fill=\extcolor] (0,2.25) rectangle (2,2.15);
          \fill[fill=\extcolor] (0,-0.25) rectangle (2,-0.15);
          \fill[fill=\intcolor] (0.1,0.1) rectangle (1.9,1.9);
      \end{scope}
      \begin{scope}[xshift=2.4cm]
          \filldraw[fill=white, draw=black, very thick, aligned dash={dash on=4pt, dash off=3pt}] (0,0) rectangle (2,2);
          \fill[fill=\extcolor] (-0.25,0) rectangle (-0.15,2);
          \fill[fill=\extcolor] (2.25,0) rectangle (2.15,2);
          \fill[fill=\extcolor] (0,2.25) rectangle (2,2.15);
          \fill[fill=\extcolor] (0,-0.25) rectangle (2,-0.15);
          \fill[fill=\intcolor] (0.1,0.1) rectangle (1.9,1.9);
      \end{scope}
    \end{tikzpicture}}}
  \caption{Two interpretations of non-overlapping domain decomposition for nodal and modal discretizations, with interface data (red) and interior data (blue). (a) In a nodal discretization, neighboring elements communicate directly through degrees of freedom at nodes, which can be partitioned into shared interface nodes and local interior nodes. (b) In a modal discretization, neighboring elements communicate indirectly through unshared interface functions, allowing for coefficients in a modal discretization to be spatially separated.}
  \label{fig:vals_vs_coeffs}
\end{figure}
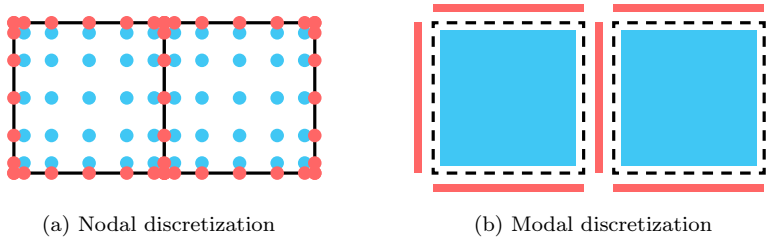

Adapting a domain decomposition approach such as the HPS scheme---originally formulated around a spectral collocation method~\cite{Martinsson_09_01, Martinsson_13_01}---to a modal discretization such as the ultraspherical spectral method gives rise to a few subtleties. In the nodal setting, values along interfaces are inherently shared between elements, allowing for an intuitive way to separate the nodes in each element into ``interior'' and ``interface'' degrees of freedom and solve for them accordingly (see \cref{fig:vals_vs_coeffs:nodal}). Cross point conditions (e.g., at a point where the corners of four quadrilaterals meet) can then be avoided by removing the degrees of freedom located at cross points~\cite{Babb_18_02}. In the modal setting, on the other hand, the coefficients in a bivariate Chebyshev expansion are not spatially localized, and therefore do not intuitively separate into such categories. To regain a decoupling for Chebyshev coefficients, it is helpful to think about bivariate functions on each element communicating not with each other directly, but with univariate functions on each interface (see \cref{fig:vals_vs_coeffs:modal}). Using a modal discretization for these bivariate interior functions and univariate interface functions then allows Chebyshev coefficients to be separated as before. Cross point conditions must then be imposed directly for the resulting linear systems to be nonsingular (see \cref{sec:ultra_hps:build}).

An alternative remedy to localize modal discretizations is to use a basis that has intrinsic spatial separation between interior and interface, such as a basis consisting of bubble functions (functions that are zero on the edges of an element) and edge functions (functions that are nonzero on the edges of an element)~\cite{Sherwin_05_01}. However, such a basis may not yield a sparse discretization of the PDE. We choose to use the ultraspherical basis to obtain sparse linear algebra, which affords our method a lower computational complexity with respect to $p$.

\subsection{Model problem: two ``glued'' squares}\label{sec:ultra_dd}
\begin{figure}[htb]
\begin{center}
\begin{tikzpicture}[scale=1.2]
\filldraw[fill=blue!10, draw=black,very thick] (0,0) rectangle (3,3);
\filldraw[fill=blue!10, draw=black,very thick] (3,0) rectangle (6,3);
\node[] at (3,3.4){\large $\Gamma$};
\node[] at (1.5,1.5){$\nabla^2 u_1 = f_1$};
\node[] at (4.5,1.5){$\nabla^2 u_2 = f_2$};
\node[] at (4.9,3.8){$\begin{aligned}&\qquad u_1 = u_2 \\ &\qquad \tfrac{\partial u_1}{\partial \vec{n}_1} + \tfrac{\partial u_2}{\partial \vec{n}_2} = 0\end{aligned}$};
\draw[thick,-latex] (4.05,3.8) to [out=200,in=30] (3.05,2);
\draw[very thick,-latex] (3,0.7) to (3.8,0.7);
\draw[very thick,-latex] (3,0.5) to (2.2,0.5);
\node[] at (4,0.6){$\vec{n}_1$};
\node[] at (2,0.4){$\vec{n}_2$};
\draw[thick,-latex] (-1,2) to [out=270,in=150] (-0.05,1.5);
\node[] at (-1,2.2){$u_1 = g_1$};
\draw[thick,-latex] (7,2) to [out=270,in=30] (6.05,1.5);
\node[] at (7,2.2){$u_2 = g_2$};
\end{tikzpicture}
\end{center}
\caption{The canonical problem setup for two ``glued'' squares.}
\label{fig:two_squares}
\end{figure}
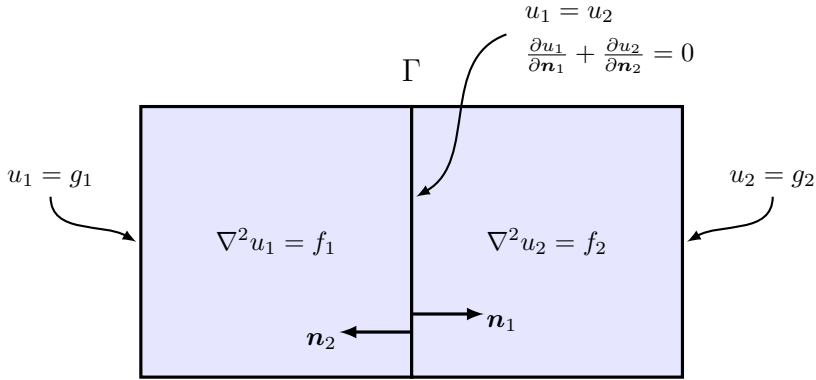

To begin, we consider the simple domain decomposition setting of two square-shaped elements that are ``glued'' together. That is, we wish to use the ultraspherical spectral method to solve the patching problem\footnote{It is worth noting that this formulation is equivalent to the global problem $\nabla^2 u = f$ in $\Omega$, $u=g$ on $\partial\Omega$, for any domain $\Omega$. This holds for any second-order linear elliptic boundary value problem~\cite{Canuto_07_01}.}
\begin{equation}\label{eq:model}
\begin{aligned}
\nabla^2 u_1 &= f_1 && \text{ in } \mathcal{E}_1, \\
\nabla^2 u_2 &= f_2 && \text{ in } \mathcal{E}_2, \\
u_1 &= g_1 && \text{ on } \partial\mathcal{E}_1 \cap \partial\Omega, \\
u_2 &= g_2 && \text{ on } \partial\mathcal{E}_2 \cap \partial\Omega, \\
u_1 &= u_2 && \text{ on } \Gamma, \\
\tfrac{\partial u_1}{\partial \vec{n}_1} &+ \tfrac{\partial u_2}{\partial \vec{n}_2} = 0 && \text{ on } \Gamma,
\end{aligned}
\end{equation}
where $\mathcal{E}$ is a mesh of the domain $\Omega = [-2, 2] \times [-1, 1]$ with elements $\mathcal{E}_1 = [-2, 0] \times [-1, 1]$ and $\mathcal{E}_2 = [0, 2] \times [-1, 1]$, $\Gamma$ is the interface between the two elements, $f$ and $g$ are given functions, and $f_i = f|_{\mathcal{E}_i}$ for any function $f$. This model problem of a pairwise merge serves as a building block in the HPS scheme. The problem setup is depicted in \cref{fig:two_squares}.

The patching problem \cref{eq:model} couples two three-sided Dirichlet problems via continuity conditions across the interface~$\Gamma$. Equivalently, \cref{eq:model} can be regarded as two decoupled, four-sided Dirichlet problems when given a suitable piece of Dirichlet data along $\Gamma$. That is, there exists an interface function $\varphi$ such that \cref{eq:model} is equivalent to
\begin{equation}\label{eq:decoupled}
\begin{aligned}
\nabla^2 u_1 &= f_1 && \text{ in } \mathcal{E}_1, &\qquad \nabla^2 u_2 &= f_2 && \text{ in } \mathcal{E}_2, \\
u_1 &= g_1 && \text{ on } \partial\mathcal{E}_1 \cap \partial\Omega, &\qquad u_2 &= g_2 && \text{ on } \partial\mathcal{E}_2 \cap \partial\Omega, \\
u_1 &= \varphi && \text{ on } \Gamma, &\qquad u_2 &= \varphi && \text{ on } \Gamma.
\end{aligned}
\end{equation}
To determine this unknown interface function $\varphi$, we aim to build a direct solver---an operator $S_\Gamma$ such that $\varphi = S_\Gamma \, g$---using ingredients from local operators on each element. In particular, we construct local direct solvers on $\mathcal{E}_1$ and $\mathcal{E}_2$, and then use pieces of these operators to construct the interfacial solution operator $S_\Gamma$. Once the interface function $\varphi$ is found, the two subproblems in \cref{eq:decoupled} decouple and can be solved independently by applying local direct solvers on $\mathcal{E}_1$ and $\mathcal{E}_2$. By building a direct solver for the global interface problem based on direct solvers for the subproblems in \cref{eq:decoupled}, the generalization to multiple elements follows naturally.

\subsubsection{Constructing local operators}\label{sec:construct_ops}

\begin{figure}[htb]
	\centering
    	\subfloat[Solution operator]{%
    		\label{fig:operators:sol}
        	\begin{minipage}[c][2.5cm][c]{0.47\textwidth}
        		\hspace{0.2cm}
                \begin{tikzpicture}[scale=0.8]
                \begin{scope}[xshift=-0.0cm]
                \filldraw[fill=white, draw=black,very thick, aligned dash={dash on=4pt, dash off=3pt}] (0,0) rectangle (2,2);
                \fill[fill=\extcolor] (-0.25,0) rectangle (-0.15,2);
                \fill[fill=\extcolor] (2.25,0) rectangle (2.15,2);
                \fill[fill=\extcolor] (0,2.25) rectangle (2,2.15);
                \fill[fill=\extcolor] (0,-0.25) rectangle (2,-0.15);
                \end{scope}
                \node[] at (3,1){\Large$\mapsto$};
                \begin{scope}[xshift=-0.2cm]
                \filldraw[fill=white, draw=black,very thick, aligned dash={dash on=4pt, dash off=3pt}] (4,0) rectangle (6,2);
                \fill[fill=\intcolor] (4.1,0.1) rectangle (5.9,1.9);
                \end{scope}
                \end{tikzpicture}\vspace{0.1cm}
    	    \end{minipage}%
        }%
        ~~~~
    	\subfloat[Dirichlet-to-Neumann operator]{%
    		\label{fig:operators:dtn}
        	\begin{minipage}[c][2.5cm][c]{0.48\textwidth}
                \hspace{0.2cm}
                \begin{tikzpicture}[scale=0.8]
                \begin{scope}[xshift=-0.0cm]
                \filldraw[fill=white, draw=black,very thick, aligned dash={dash on=4pt, dash off=3pt}] (0,0) rectangle (2,2);
                \fill[fill=\extcolor] (-0.25,0) rectangle (-0.15,2);
                \fill[fill=\extcolor] (2.25,0) rectangle (2.15,2);
                \fill[fill=\extcolor] (0,2.25) rectangle (2,2.15);
                \fill[fill=\extcolor] (0,-0.25) rectangle (2,-0.15);
                \end{scope}
                \node[] at (3,1){\Large$\mapsto$};
                \begin{scope}[xshift=0.3cm]
                \filldraw[fill=white, draw=black,very thick, aligned dash={dash on=4pt, dash off=3pt}] (4,0) rectangle (6,2);
                \fill[fill=\extcolor] (3.75,0) rectangle (3.85,2);
                \fill[fill=\extcolor] (6.25,0) rectangle (6.15,2);
                \fill[fill=\extcolor] (4,2.25) rectangle (6,2.15);
                \fill[fill=\extcolor] (4,-0.25) rectangle (6,-0.15);
                \draw[\extcolor,very thick,-latex] (6.25,1) to (6.75,1);
                \draw[\extcolor,very thick,-latex] (6.25,1.5) to (6.75,1.5);
                \draw[\extcolor,very thick,-latex] (6.25,0.5) to (6.75,0.5);
                \draw[\extcolor,very thick,-latex] (3.75,1) to (3.25,1);
                \draw[\extcolor,very thick,-latex] (3.75,1.5) to (3.25,1.5);
                \draw[\extcolor,very thick,-latex] (3.75,0.5) to (3.25,0.5);
                \draw[\extcolor,very thick,-latex] (5,2.25) to (5,2.75);
                \draw[\extcolor,very thick,-latex] (4.5,2.25) to (4.5,2.75);
                \draw[\extcolor,very thick,-latex] (5.5,2.25) to (5.5,2.75);
                \draw[\extcolor,very thick,-latex] (5,-0.25) to (5,-0.75);
                \draw[\extcolor,very thick,-latex] (4.5,-0.25) to (4.5,-0.75);
                \draw[\extcolor,very thick,-latex] (5.5,-0.25) to (5.5,-0.75);
                \end{scope}
                \end{tikzpicture}\vspace{0.1cm}
    	    \end{minipage}%
        }%
	\caption{A visualization of the local operators computed for each element. (a) The solution operator on a quadrilateral takes in four univariate functions of Dirichlet data and returns a bivariate function that solves the PDE using the given boundary conditions. (b) The Dirichlet-to-Neumann operator on a quadrilateral takes in four univariate functions of Dirichlet data and returns four univariate functions of Neumann data, representing the normal derivative of the solution to the PDE on the four sides.}
	\label{fig:operators}
\end{figure}
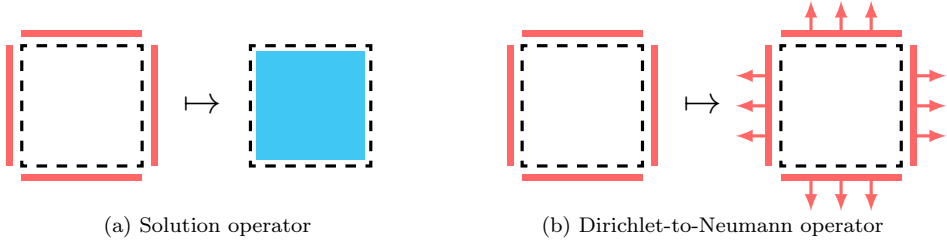

To construct a direct solver for \cref{eq:decoupled}, we first build operators that encode how to solve the PDE locally on elements $\mathcal{E}_1$ and $\mathcal{E}_2$. Such operators, called solution operators, take in Dirichlet data and return the corresponding solution to the PDE on an element. For a quadrilateral domain, the solution operator takes in four univariate functions---representing four sides of Dirichlet data---and returns a bivariate function that satisfies the PDE (see \cref{fig:operators:sol}).

We use the ultraspherical spectral method for solving PDEs on quadrilaterals (see \cref{sec:quad_tri}). If on each element we employ a $(p+1) \times (p+1)$ coefficient discretization for the solution so that the solution is at most a degree-$(p,p)$ polynomial, then the solution operator $S_\mathcal{E} \in \mathbb{C}^{(p+1)^2 \times 4(p+1)+1}$ on element $\mathcal{E}$ is a dense matrix. For a column vector $\vec{c} \in \mathbb{C}^{4(p+1)}$ and scalar $\alpha \in \mathbb{C}$, the product $S_\mathcal{E} \,\begin{bsmallmatrix} \vec{c} \\ \alpha \end{bsmallmatrix} \in \mathbb{C}^{(p+1)^2}$ represents the $(p+1) \times (p+1)$ Chebyshev coefficients of the solution to the PDE on $\mathcal{E}$ with Dirichlet data $\vec{c}$ and righthand side $\alpha f$.\footnote{In practice, we always take $\alpha = 1$.} Here, $\vec{c} = [\vec{c}_1, \vec{c}_2, \vec{c}_3, \vec{c}_4]^T$ represents the Chebyshev coefficients of four univariate functions of Dirichlet data on the left, right, bottom, and top of $\mathcal{E}$, respectively, each discretized with $p+1$ coefficients. For example, on the left side, the coefficients $\vec{c}_1 \in \mathbb{C}^{p+1}$ define the degree-$p$ boundary function $h_1(y)$ as
\[
h_1(y) = \sum_{j=0}^p (\vec{c}_1)_j T_j(y).
\]
Similarly, $\vec{c}_2$, $\vec{c}_3$, and $\vec{c}_4$ define functions on the other three sides.

The solution operator on $\mathcal{E}$ can be decomposed into four operators $S_\mathcal{E}^1, S_\mathcal{E}^2, S_\mathcal{E}^3, S_\mathcal{E}^4 \in \mathbb{C}^{(p+1)^2 \times (p+1)}$ that account for the homogeneous part of the solution and one column vector $S_\mathcal{E}^\text{rhs} \in \mathbb{C}^{(p+1)^2 \times 1}$ that accounts for the particular part of the solution.\footnote{Although including the particular solution in the solution operator is not typically standard in the HPS literature, we do this here because it avoids repeating a description of the linear algebra when constructing the particular solution separately and matches how \ultraSEM is implemented.} That is,
\[
S_\mathcal{E} = \begin{bmatrix}\!\begin{array}{cccc|c}
&&&\\[-0.3em]
S_\mathcal{E}^1 & S_\mathcal{E}^2 & S_\mathcal{E}^3 & S_\mathcal{E}^4 & S_\mathcal{E}^\text{rhs}\\[0.8em]
\end{array}\!\end{bmatrix}\!,
\]
where the vector $S_\mathcal{E}^{\text{rhs}}$ is defined by
\[
S_\mathcal{E}^\text{rhs} = \operatorname{vec}(X), \qquad u_\text{rhs}(x,y) = \sum_{i=0}^p \sum_{j=0}^p X_{ij} T_i(y) T_j(x),
\]
and $\operatorname{vec}(\cdot)$ is the column-wise vectorization operator. Here, $u_\text{rhs}$ satisfies the PDE on $\mathcal{E}$ with homogeneous boundary conditions, i.e.,
\[
\nabla^2 u_\text{rhs} = f|_\mathcal{E}, \qquad u_\text{rhs}|_{\partial\mathcal{E}} = 0.
\]
The products $S_\mathcal{E}^i \vec{c}_i$ represent the $(p+1) \times (p+1)$ Chebyshev coefficients of the approximate solution to the homogeneous problem (i.e., $f = 0$) with Dirichlet data on side $i$ given by the Chebyshev coefficients $\vec{c}_i$ and zero Dirichlet data on the other three sides.
Thus, the solution operator $S_\mathcal{E}$ depends on the domain $\mathcal{E}$, the PDO, and the righthand side $f$, but not the Dirichlet data $g$. However, the solution operator can be efficiently updated if $f$ is changed (see the discussion of \mytt{updateRHS} in \cref{sec:software}).

We construct the matrices $S_\mathcal{E}^i$ column-by-column. To construct the $j$th column of $S_\mathcal{E}^i$, we set the $j$th Dirichlet coefficient to one and the rest to zero, i.e.,
\begin{equation}\label{eq:delta_coeffs}
(\vec{c}_i)_k = \begin{cases}
1, & \text{if } k=j, \\
0, & \text{otherwise,}
\end{cases}
\end{equation}
for $0 \leq k \leq p$. We wish to solve the PDE using this Dirichlet data for each $0 \leq j \leq p$ to obtain the $(p+1) \times (p+1)$ coefficients of the solution, which are reshaped and placed as a column into $S_\mathcal{E}^i$. That is, the $j$th column of the solution operator for the $i$th side of the element $\mathcal{E}$, i.e., $(S_\mathcal{E}^i)_{:,j}$, is constructed as
\[
(S_\mathcal{E}^i)_{:,j} = \operatorname{vec}(X), \qquad v_j(x,y) = \sum_{k=0}^p \sum_{\ell=0}^p X_{k\ell} T_k(y) T_\ell(x),
\]
where $v_j$ approximately satisfies the following homogeneous PDE:
\[
\nabla^2 v_j = 0, \qquad v_j|_{\partial\mathcal{E}} = \begin{cases} T_j, & \text{on side } i, \\ 0, & \text{otherwise}. \end{cases}
\]
Unfortunately, the Dirichlet data used in this construction process may have discontinuities at the corners of the domain, leading to incompatible boundary conditions. To ensure compatibility is satisfied, we orthogonally project each function $\vec{c}_i$ onto the space of functions that are continuous at the corners before solving the PDE. The compatibility conditions at the four corners of the quadrilateral can be encoded into a matrix $B \in \mathbb{C}^{4 \times 4(p+1)}$ given by
\[
\hspace{1.6cm}
B = \begin{bmatrix}
B_{-1} & 0 & -B_{-1} & 0 \\
B_{+1} & 0 & 0 & -B_{-1} \\
0 & B_{-1} & -B_{+1} & 0 \\
0 & B_{+1} & 0 & -B_{+1}
\end{bmatrix}%
\!\!\!\!
\begin{array}{l}
\}\,\text{\footnotesize bottom left corner} \\
\}\,\text{\footnotesize top left corner} \\
\}\,\text{\footnotesize bottom right corner} \\
\}\,\text{\footnotesize top right corner}
\end{array}
\]
with
\[
B_{\pm 1} =
\begin{bmatrix}
T_0(\pm1) & T_1(\pm1) & \cdots & T_p(\pm1)
\end{bmatrix},
\]
where $T_j(\pm1) = (\pm1)^j$. The matrix $B_{\pm 1}$ is an evaluation operator at the endpoints of the interval $[-1, 1]$. So, for the functions $h_i$ defined above, $B_{\pm 1} \vec{c}_i = h_i(\pm 1)$. A given piece of boundary data defined by the coefficients $\vec{c} = [\vec{c}_1, \vec{c}_2, \vec{c}_3, \vec{c}_4]$ is compatible at the corners if and only if $B\vec{c} = 0$. To project the boundary data so that it satisfies compatibility, we build a basis for $\operatorname{null}(B)$, which is of rank $4(p+1)-4$. Taking the singular value decomposition $B = U\Sigma V^*$ and letting $\tilde{V}$ be the last $4(p+1)-4$ columns of $V$, we construct a projection matrix $P = \tilde{V} \tilde{V}^*$. Since this projection matrix depends only on $p$, it can be precomputed and stored. The product $\tilde{\vec{c}} = P\vec{c}$ orthogonally projects the functions defined by $\vec{c}_1$, $\vec{c}_2$, $\vec{c}_3$, and $\vec{c}_4$ onto the space of compatible boundary conditions, so that $\tilde{\vec{c}}_1$, $\tilde{\vec{c}}_2$, $\tilde{\vec{c}}_3$, and $\tilde{\vec{c}}_4$ are continuous at the four corners of the quadrilateral. We apply this projection during the construction process to the Dirichlet data $\vec{c}$ in \cref{eq:delta_coeffs} to obtain compatible Dirichlet data $\tilde{\vec{c}}$. It is this Dirichlet data that we use to construct the columns of the solution operator $S_{\mathcal{E}}^i$.

Continuity conditions between elements are communicated locally via the Dirichlet-to-Neumann operator, or Poincar\'{e}--Steklov operator. The Dirichlet-to-Neumann operator on an element $\mathcal{E}$, denoted by $\dtn_\mathcal{E}$, maps Dirichlet data on each side of $\mathcal{E}$ to the outward flux of the local solution to the PDE on each side of $\mathcal{E}$ (see \cref{fig:operators:dtn}). One may apply $\dtn_\mathcal{E}$ by first computing the local solution to the PDE on $\mathcal{E}$ for the given Dirichlet data and then evaluating the outward flux of the solution on the boundary. Hence, the Dirichlet-to-Neumann operator can be written as a product of the normal derivative operator and the solution operator. That is, $\dtn_\mathcal{E} = D_\mathcal{E} S_\mathcal{E}$, where $D_\mathcal{E}$ computes the outward flux of a bivariate function on each side of the element $\mathcal{E}$ when given its $(p+1)^2$ Chebyshev coefficients. On the reference square $[-1,1]^2$, $D_{[-1,1]^2} \in \mathbb{C}^{4(p+1) \times (p+1)^2}$ is given by
\[
D_{[-1,1]^2}
=
\begin{bmatrix}
I \otimes D_{-1} \\
I \otimes D_{+1} \\
D_{-1} \otimes I \\
D_{+1} \otimes I
\end{bmatrix}%
\!\!\!\!
\begin{array}{l}
\}\,\text{\footnotesize left normal derivative} \\
\}\,\text{\footnotesize right normal derivative} \\
\}\,\text{\footnotesize bottom normal derivative} \\
\}\,\text{\footnotesize top normal derivative}
\end{array}
\]
where `$\otimes$' denotes the Kronecker product operator for matrices, $I$ is the $(p+1) \times (p+1)$ identity matrix, and
\[
D_{\pm 1}
= \pm
\begin{bmatrix}
T'_0(\pm1) & T'_1(\pm1) & \cdots & T'_p(\pm1)\end{bmatrix}, \qquad T'_j(\pm1)= (\pm1)^j j^2.
\]
On quadrilaterals and triangles, the normal derivative operator is transformed according to the Jacobian factors described in \cref{sec:quad_tri}. Hence, the Dirichlet-to-Neumann operator $\dtn_\mathcal{E} \in \mathbb{C}^{4(p+1) \times (4(p+1)+1)}$ is a dense matrix. The product $\dtn_\mathcal{E} \,\begin{bsmallmatrix} \vec{c} \\ \alpha \end{bsmallmatrix} \in \mathbb{C}^{4(p+1)}$ represents the four normal derivatives of the solution to the PDE on the element $\mathcal{E}$ with Dirichlet data $\vec{c}$ and righthand side $\alpha f$,\footnote{Again, note that our definition of the Dirichlet-to-Neumann operator includes the particular solution.}  each discretized with $p+1$ Chebyshev coefficients. In the context of the model problem \cref{eq:decoupled}, the Dirichlet-to-Neumann operators $\dtn_{\mathcal{E}_1}$ and $\dtn_{\mathcal{E}_2}$ on the elements $\mathcal{E}_1$ and $\mathcal{E}_2$, respectively, are merged to make the interfacial solution operator $S_\Gamma$, allowing for the direct solution of the unknown interface function $\varphi$.

\subsubsection{Merging two operators}
With local operators constructed on each element $\mathcal{E}_1$ and $\mathcal{E}_2$, we now aim to build a global solution operator, $S_\Gamma$, from the local operators $S_{\mathcal{E}_1}$, $S_{\mathcal{E}_2}$, $\dtn_{\mathcal{E}_1}$, and $\dtn_{\mathcal{E}_2}$, to solve for the unknown interface function $\varphi$. Mathematically, this decomposition mimics the classical Schur complement method for domain decomposition, keeping the physical interpretation for modal discretizations from \cref{sec:modal_dd} in mind.

For elements $\mathcal{E}_1$ and $\mathcal{E}_2$, let $\Gamma_1$ and $\Gamma_2$ denote the indices of the local Dirichlet data corresponding to the shared boundary $\Gamma$. For $\mathcal{E}_1$, since the shared interface $\Gamma$ is on the right side and the boundary data $\vec{c} = [\vec{c}_1, \vec{c}_2, \vec{c}_3, \vec{c}_4]$ is ordered as left, right, bottom, and top, the indices corresponding to the $p+1$ Chebyshev coefficients of the right-side Dirichlet data $\vec{c}_2$ are given by the set $\Gamma_1 = \{(p+1)+1, \ldots, 2(p+1)\}$. Similarly, since the interface $\Gamma$ is on the left side of $\mathcal{E}_2$, the indices of the local Dirichlet data on the shared boundary of $\mathcal{E}_2$ are given by $\Gamma_2 = \{1, \ldots, p+1\}$. Finally, denote by $\mathrm{L}_1$ and $\mathrm{L}_2$ the sets containing the indices corresponding to the coefficients of the unshared Dirichlet data on each element, so that $\mathrm{L}_1 = \{1, \ldots, 4(p+1)\} \setminus \Gamma_1$ and $\mathrm{L}_2 = \{1, \ldots, 4(p+1)\} \setminus \Gamma_2$. For a matrix

With these indices defined for $\mathcal{E}_1$ and $\mathcal{E}_2$ based on interaction with the Dirichlet data on $\Gamma$, the rows and columns of the local operators $\dtn_{\mathcal{E}_1}$ and $\dtn_{\mathcal{E}_2}$ can be partitioned into ``interior'' and ``interface'' blocks. The pieces of $\dtn_{\mathcal{E}_1}$ and $\dtn_{\mathcal{E}_2}$ that affect the shared interface naturally separate, and a Schur complement may be performed to write down the following  $(p+1) \times (p+1)$ linear system for the solution operator on the interface:
\begin{equation}\label{eq:sol_2}
-\left(\dtn^{\Gamma_1,\Gamma_1}_{\mathcal{E}_1} + \dtn^{\Gamma_2,\Gamma_2}_{\mathcal{E}_2}\right) S_\Gamma = \left[\!\begin{array}{>{\centering\arraybackslash$}m{1cm}<{$}>{\centering\arraybackslash$}m{1cm}<{$}|>{\centering\arraybackslash$}m{2.7cm}<{$}} && \\[-1em] \dtn^{{\Gamma_1,\mathrm{L}_1}}_{{\mathcal{E}_1}} & \!\dtn^{\Gamma_2,\mathrm{L}_2}_{\mathcal{E}_2} & \;\dtn^{\Gamma_1,\text{end}}_{\mathcal{E}_1} + \dtn^{\Gamma_2,\text{end}}_{\mathcal{E}_2} \\[-1em] && \end{array}\!\right],
\end{equation}
where the last column of the righthand side of \cref{eq:sol_2} encodes the contribution from the particular solution. Here, superscripts denote row and column indices for slicing a matrix and ``end'' denotes the index of the last column of a matrix. The linear system in \cref{eq:sol_2} has a clear interpretation: the matrix $\dtn^{\Gamma_1,\Gamma_1}_{\mathcal{E}_1} + \dtn^{\Gamma_2,\Gamma_2}_{\mathcal{E}_2}$ computes the jump in the normal derivative across the shared interface $\Gamma$ and enforces this jump to be offset by the contributions from the unshared sides and particular solution, resulting in a discrete analogue of the original continuity condition in \cref{eq:model}. As before, the merged solution operator $S_\Gamma \in \mathbb{C}^{(p+1) \times (6(p+1)+1)}$ is a dense matrix. For a column vector $\vec{c} \in \mathbb{C}^{6(p+1)}$ and scalar $\alpha \in \mathbb{C}$, the product $S_\Gamma \,\begin{bsmallmatrix} \vec{c} \\ \alpha \end{bsmallmatrix} \in \mathbb{C}^{p+1}$ represents the $p+1$ Chebyshev coefficients of the solution to the PDE on $\Gamma$ with Dirichlet data $\vec{c} = [\vec{c}_1, \ldots, \vec{c}_6]^T$ and righthand side $\alpha f$, where now the Dirichlet data $\vec{c}$ is specified on the six sides of the merged domain $\Omega$.

The Schur complement also allows us to write down the Dirichlet-to-Neumann operator for the merged domain. Using the new solution operator $S_\Gamma$, we can construct a new Dirichlet-to-Neumann operator on $\Omega$ as
\begin{equation}\label{eq:dtn_2}
\dtn_\Omega = \left[\!\begin{array}{cc|c}
        && \\[-1em]
		\dtn^{\mathrm{L}_1,\mathrm{L}_1}_{\mathcal{E}_1} & 0 & \,\dtn^{\mathrm{L}_1,\text{end}}_{\mathcal{E}_1} \\[0.5em]
    	0 & \dtn^{\mathrm{L}_2,\mathrm{L}_2}_{\mathcal{E}_2} & \,\dtn^{\mathrm{L}_2,\text{end}}_{\mathcal{E}_2} \\[-1em] &&
	\end{array}\!\right] + \left[\!\begin{array}{c} \\[-1em] \dtn^{\mathrm{L}_1,\Gamma_1}_{\mathcal{E}_1} \\[0.4em] \dtn^{\mathrm{L}_2,\Gamma_2}_{\mathcal{E}_2} \\[0.4em] \end{array}\!\right] S_\Gamma,
\end{equation}
where $\dtn_\Omega \in \mathbb{C}^{6(p+1) \times (6(p+1)+1)}$. The vector $\dtn_\Omega \,\begin{bsmallmatrix} \vec{c} \\ \alpha \end{bsmallmatrix}$ represents normal derivatives on the six sides of $\Omega$ of the solution to the PDE on $\Omega$ with Dirichlet data $\vec{c}$ and righthand side $\alpha f$, each discretized with $p+1$ Chebyshev coefficients.

\subsubsection{Computing the solution}

We now have all the ingredients we need to compute the solution to \cref{eq:model}. We begin by converting the given boundary functions, $g_1$ and $g_2$, into Chebyshev coefficients. On each of the three sides of $\mathcal{E}_1$ and $\mathcal{E}_2$ where $g_1$ and $g_2$ are known, we construct the degree-$p$ Chebyshev approximant to the boundary data and compile the coefficients into vectors $\vec{g}_1$ and $\vec{g}_2$ of length $3(p+1)$. Next, to solve for the interface function $\varphi$ that makes \cref{eq:decoupled} equivalent to \cref{eq:model}, we simply compute the matrix-vector product $S_\Gamma \, \begin{bsmallmatrix} \vec{g} \\ 1 \end{bsmallmatrix}$, where $\vec{g} = [\vec{g}_1, \vec{g}_2]^T$, which yields the $p+1$ Chebyshev coefficients of $\varphi$. With the Dirichlet data now known on all four sides of each of the elements $\mathcal{E}_1$ and $\mathcal{E}_2$, the local solution operators $S_{\mathcal{E}_1}$ and $S_{\mathcal{E}_2}$ can finally be applied. Defining vectors $\widehat{\vec{g}}_i$ such that $\widehat{\vec{g}}_i^{\,\Gamma_i} = \varphi$ and $\widehat{\vec{g}}_i^{\,\mathrm{L}_i} = \vec{g}_i$ for $i=1, 2$, the matrix-vector products $S_{\mathcal{E}_1} \,\begin{bsmallmatrix} \widehat{\vec{g}}_1 \\ 1 \end{bsmallmatrix}$ and $S_{\mathcal{E}_2} \,\begin{bsmallmatrix} \widehat{\vec{g}}_2 \\ 1 \end{bsmallmatrix}$ contain the $(p+1) \times (p+1)$ coefficients of the solutions $u_1$ and $u_2$, respectively, satisfying \cref{eq:model}.

\subsection{The hierarchical scheme}\label{sec:ultra_hps}

At the end of merge process for the model problem of two ``glued'' squares, we are left with two operators acting on $\Omega$: (1) a solution operator, $S_\Gamma$, to solve for the unknown interface inside $\Omega$, and (2) a Dirichlet-to-Neumann operator, $\dtn_\Omega$, to map boundary data to outward fluxes on $\Omega$. These operators encode everything we need to know to solve the PDE on $\Omega$. In effect, $\Omega$ is now no different from the original elements $\mathcal{E}_1$ or $\mathcal{E}_2$, and so it can be treated as just another element, ready to be merged again with a new domain. After another merge, we are once again in the same situation, with access to local operators that allow us to treat the merged domain as a black box. This is the hierarchical Poincar\'{e}--Steklov scheme.

\subsubsection{Build stage}\label{sec:ultra_hps:build}

For a mesh $\mathcal{E} = \{\mathcal{E}_i\}_{i=1}^{n_\text{elem}}$ of a domain $\Omega$, the scheme begins with a local build stage, wherein local solution operators $S_{\mathcal{E}_i}$ and Dirichlet-to-Neumann operators $\dtn_{\mathcal{E}_i}$ are constructed on each element $\mathcal{E}_i$ according to \cref{sec:construct_ops}. The local build process is outlined in \cref{alg:init} and is referred to as ``initialization'' in \ultraSEM. As the operations performed in this stage are local to each element, \cref{alg:init} can be parallelized across elements.

\begin{algorithm}[htb]
\renewcommand\thealgorithm{\thesection.1\alph{algorithm}}
\caption{Local build stage: $\texttt{initialize}(\mathcal{E}, \mathcal{L}, f)$}
\begin{algorithmic}[1]
\Require{Mesh $\mathcal{E} = \{\mathcal{E}_i\}_{i=1}^{n_\text{elem}}$, partial differential operator $\mathcal{L}$, righthand side $f$}
\Ensure{Solution and Dirichlet-to-Neumann operators for every element, $\{S_{\mathcal{E}_i},\dtn_{\mathcal{E}_i}\}_{i=1}^{n_\text{elem}}$}
\For{each element $\mathcal{E}_i$ in mesh}
    \State{Transform $(\mathcal{L}, f) \mapsto (\widehat{\mathcal{L}}, \widehat{f})$ into reference space (see \cref{sec:quad_tri}).}
	\State{Discretize $\widehat{\mathcal{L}}$ and $\widehat{f}$ using the ultraspherical spectral method.}
	\State{Construct the solution operator $S_{\mathcal{E}_i}$ (see \cref{sec:construct_ops}).}
	\State{\multiline{Construct the Dirichlet-to-Neumann operator $\dtn_{\mathcal{E}_i} := D_{\mathcal{E}_i}S_{\mathcal{E}_i}\!$ \\ (see \cref{sec:construct_ops}).}}
\EndFor
\State{\Return $\{S_{\mathcal{E}_i}, \dtn_{\mathcal{E}_i}\}_{i=1}^{n_\text{elem}}$}
\end{algorithmic}
\label{alg:init}
\end{algorithm}

Once local operators have been computed for each element, the scheme enters the global build stage, where a hierarchy of merged operators is constructed in an upward pass. Given a set of indices $\mathcal{I}$ that define a sequence of pairwise merges between elements, operators are merged as in \cref{eq:sol_2,eq:dtn_2} in the order $\mathcal{I}$ until the entire mesh has been merged into one large conglomerate. Along the way, merged elements store their newly computed solution operators and Dirichlet-to-Neumann operators. The global build stage ends with a solution operator that acts on the entire mesh, taking in Dirichlet data on every boundary of $\Omega$ and returning the solution to the PDE along the penultimate merged interface. The global build stage is outlined in \cref{alg:build}.

When the mesh contains cross points (i.e., points in the interior of the mesh where corners of multiple elements meet), the linear system defining the solution operator,
\begin{equation}\label{eq:build_sol}
-\left(\dtn^{\Gamma_i,\Gamma_i}_{\mathcal{E}_i} + \dtn^{\Gamma_j,\Gamma_j}_{\mathcal{E}_j}\right) S_{\Gamma_{ij}} = \left[\!\begin{array}{>{\centering\arraybackslash$}m{1cm}<{$}>{\centering\arraybackslash$}m{1cm}<{$}|>{\centering\arraybackslash$}m{2.7cm}<{$}} && \\[-1em] \dtn^{{\Gamma_i,\mathrm{L}_i}}_{{\mathcal{E}_i}} & \!\dtn^{\Gamma_j,\mathrm{L}_j}_{\mathcal{E}_j} & \;\dtn^{\Gamma_i,\text{end}}_{\mathcal{E}_i} + \dtn^{\Gamma_j,\text{end}}_{\mathcal{E}_j} \\[-1em] && \end{array}\!\right],
\end{equation}
may be rank deficient, as a continuity condition on the sum of the normal fluxes around the cross point has not been imposed~\cite{Canuto_07_01}. 
Rather than imposing this condition directly, we solve the rank-deficient system by projecting out the cross-point modes, which has the same effect as removing the degrees of freedom located at cross points in a collocation-based scheme~\cite{Babb_18_02}. The nullspace of \cref{eq:build_sol} contains precisely the cross-point modes, and so we implement this projection step by performing a minimum-norm least-squares solve on \cref{eq:build_sol}. As this is a projection method, the resulting residual is guaranteed to be identically zero.


\begin{algorithm}[htb]
\renewcommand\thealgorithm{\thesection.1\alph{algorithm}}
\caption{Global build stage (upward pass): $\texttt{build}(\mathcal{E}, \dtn_\mathcal{E}, \mathcal{I})$}
\begin{algorithmic}[1]
\Require{Mesh $\mathcal{E} = \{\mathcal{E}_i\}_{i=1}^{n_\text{elem}}$, local operators $\dtn_\mathcal{E} = \{\dtn_{\mathcal{E}_i}\}_{i=1}^{n_\text{elem}}$, merge indices $\mathcal{I}$}
\Ensure{Solution and Dirichlet-to-Neumann operators for every merge, $\{S_{\Gamma_{ij}},\dtn_{\mathcal{E}_{ij}}\}_{(i,j) \in \mathcal{I}}$}
\For{each pair in $(i,j) \in \mathcal{I}$}
	\State{Define the merged domain $\mathcal{E}_{ij} := \mathcal{E}_i \cup \mathcal{E}_j$.}
	\State{Define the shared interface $\Gamma_{ij} := \mathcal{E}_i \cap \mathcal{E}_j$.}
	\State{Define indices $\Gamma_i$, $\Gamma_j$ for the shared boundary $\Gamma_{ij}$ on $\mathcal{E}_i$, $\mathcal{E}_j$.}
	\State{\multiline{Define indices $\mathrm{L}_i := \overline{\Gamma_i}$, $\mathrm{L}_j := \overline{\Gamma_j}$ for the unshared boundaries on $\mathcal{E}_i$, $\mathcal{E}_j$.}}
	\State{\multiline{Solve the linear system \vspace{-0.4em}
	\[
	-\left(\dtn^{\Gamma_i,\Gamma_i}_{\mathcal{E}_i} + \dtn^{\Gamma_j,\Gamma_j}_{\mathcal{E}_j}\right) S_{\Gamma_{ij}} = \left[\!\begin{array}{>{\centering\arraybackslash$}m{1cm}<{$}>{\centering\arraybackslash$}m{1cm}<{$}|>{\centering\arraybackslash$}m{2.7cm}<{$}} && \\[-1em] \dtn^{{\Gamma_i,\mathrm{L}_i}}_{{\mathcal{E}_i}} & \!\dtn^{\Gamma_j,\mathrm{L}_j}_{\mathcal{E}_j} & \;\dtn^{\Gamma_i,\text{end}}_{\mathcal{E}_i} + \dtn^{\Gamma_j,\text{end}}_{\mathcal{E}_j} \\[-1em] && \end{array}\!\right]\vspace{-0.4em}
	\]
	for the merged solution operator $S_{\Gamma_{ij}}$.}}
	\State{\multiline{Define the merged Dirichlet-to-Neumann operator, \vspace{-0.4em}
	  \begin{equation*}\label{eq:build_dtn}
	  \dtn_{\mathcal{E}_{ij}} := \left[\!\begin{array}{cc|c}
        && \\[-1em]
		\dtn^{\mathrm{L}_i,\mathrm{L}_i}_{\mathcal{E}_i} & 0 & \,\dtn^{\mathrm{L}_i,\text{end}}_{\mathcal{E}_i} \\[0.5em]
    	0 & \dtn^{\mathrm{L}_j,\mathrm{L}_j}_{\mathcal{E}_j} & \,\dtn^{\mathrm{L}_j,\text{end}}_{\mathcal{E}_j} \\[-1em] &&
	\end{array}\!\right] + \left[\!\begin{array}{c} \\[-1em] \dtn^{\mathrm{L}_i,\Gamma_i}_{\mathcal{E}_i} \\[0.4em] \dtn^{\mathrm{L}_j,\Gamma_j}_{\mathcal{E}_j} \\[0.4em] \end{array}\!\right] S_{\Gamma_{ij}}.
	\end{equation*}}}
\EndFor
\vspace{-1em}
\State{\Return $\{S_{\Gamma_{ij}}, \dtn_{\mathcal{E}_{ij}}\}_{(i,j) \in \mathcal{I}}$}
\end{algorithmic}
\label{alg:build}
\end{algorithm}

\subsubsection{Solve stage}\label{sec:ultra_hps:solve}
The final stage of the scheme is the solve stage, which uses the merged solution operators to recover the unknown interface data in a downward pass through the hierarchy. Beginning at the top of the hierarchy, the solution operator acting on the entire mesh is applied to the known Dirichlet data $g$, returning the Chebyshev coefficients of the solution on the top-level merged interface. These coefficients are then used as Dirichlet data on the next level, where solution operators are again applied to compute the unknown interface data on subdomains. Finally, at the bottom level of the hierarchy---where the solution is now known at each interface between elements---the local solution operators $S_{\mathcal{E}_i}$ are applied to compute the bivariate solution in the interior of each element $\mathcal{E}_i$. The solve stage is outlined in \cref{alg:solve}.

The solve stage may be executed multiple times using different boundary data without recomputing the operators constructed in the build stage. The stored operators may also be efficiently updated to solve \cref{eq:pde} with a different righthand side $f$. Recall that the last column of every solution operator and Dirichlet-to-Neumann operator in the hierarchy corresponds to the contribution from the particular solution. Using a new righthand side, an updated particular solution can be constructed on each element $\mathcal{E}_i$ as in \cref{sec:construct_ops}, replacing the last columns of $S_{\mathcal{E}_i}$ and $\dtn_{\mathcal{E}_i}$. A modified build stage may be then be executed, where the last column of each interfacial solution and Dirichlet-to-Neumann operator is updated by solving the linear system \cref{eq:build_sol} in an upward pass.

\begin{algorithm}[htb]
\renewcommand\thealgorithm{\thesection.2}
\caption{Solve stage (downward pass): $\texttt{solve}(\mathcal{E}, g)$}
\begin{algorithmic}[1]
\Require{Element (or merged element) $\mathcal{E}$, Dirichlet data $g$}
\Ensure{Solutions $\{u_i\}_{i=1}^{n_\text{elem}}$}
\If{$\mathcal{E}$ is the entire domain}
	\State{Get all boundary faces $(\partial \mathcal{E})_i$.}
	\State{Evaluate Dirichlet data $g_i := g((\partial \mathcal{E})_i)$ and convert to Chebyshev coefficients.}
\EndIf
\If{$\mathcal{E}$ is a leaf}
	\State{Compute the local solution $u := S_{\mathcal{E}} \begin{bsmallmatrix} g \\ 1 \end{bsmallmatrix}$.}
	\State{\Return $u$}
\Else
	\State{Look up the elements $\mathcal{E}_i$, $\mathcal{E}_j$ that were merged to make $\mathcal{E}$.}
	\State{Define the shared interface $\Gamma_{ij} := \mathcal{E}_i \cap \mathcal{E}_j$.}
	\State{Recover the missing interface data $\varphi := S_{\Gamma_{ij}} \,\begin{bsmallmatrix} g \\ 1 \end{bsmallmatrix}$.}
	\State{Define vectors $\widehat{\vec{g}}_i$, $\widehat{\vec{g}}_j$ such that $\widehat{\vec{g}}_i^{\,\Gamma_i} = \varphi$, $\widehat{\vec{g}}_i^{\,\mathrm{L}_i} = g_i$ and $\widehat{\vec{g}}_j^{\,\Gamma_j} = \varphi$, $\widehat{\vec{g}}_j^{\,\mathrm{L}_j} = g_j$.}
	\State{Compute the solution on $\mathcal{E}_i$, $\{u_i\} := \texttt{solve}(\mathcal{E}_i, \widehat{\vec{g}}_i)$.}
	\State{Compute the solution on $\mathcal{E}_j$, $\{u_j\} := \texttt{solve}(\mathcal{E}_j, \widehat{\vec{g}}_j)$.}
	\State{\Return $\{u_i\} \cup \{u_j\}$}
\EndIf
\end{algorithmic}
\label{alg:solve}
\end{algorithm}

\subsection{Computational complexity}\label{sec:complexity}
We now determine the computational complexity of the build and solve stages in terms of the number of degrees of freedom, $N \approx (p/h)^2$, where $h$ is the minimum mesh size and $p$ is the polynomial order. Here, we assume that the number of elements in the mesh, $n_\text{elem}$, scales as $\mathcal{O}(1/h^2)$, which is valid for a mesh that is approximately uniformly refined. For a mesh that is adaptively refined, the number of elements is typically much less than this estimate.

We begin with the local build stage. On each element $\mathcal{E}_i$, we approximate the solution as a degree-$(p,p)$ polynomial using $(p+1)^2$ degrees of freedom. After transforming the PDE into the local coordinate system of the element, we discretize $\widehat{\mathcal{L}}$ and $\widehat{f}$ using the ultraspherical spectral method. The bivariate Chebyshev coefficients of $\widehat{f}$ can be computed in $\mathcal{O}(p^2 \log p)$ operations via a discrete cosine transform~\cite{ATAP}. A separable representation of $\widehat{\mathcal{L}}$ can be computed in $\mathcal{O}(p^3)$ operations using the singular value decomposition, and differentiation, conversion, and multiplication matrices can be constructed for each separable piece in $\mathcal{O}(p)$ operations. The $(p+1)^2 \times (p+1)^2$ discrete PDO can then be assembled using Kronecker products in $\mathcal{O}(p^4)$ operations. The discrete PDO $L$ is almost block-banded, with a bandwidth\footnote{The bandwidth of the discrete PDO depends on the polynomial degree, $m$, used to approximate the variable coefficients. As in \cref{sec:ultraS}, we assume that $m \ll p$ so that the discrete PDO is sparse.} of $\mathcal{O}(p)$ and $\mathcal{O}(p)$ dense rows. To compute the solution operator $S_{\mathcal{E}_i}$, we must solve a linear system with $\mathcal{O}(p)$ righthand sides. That is, we must solve a system of the form $LX=B$, where $L$ is $\mathcal{O}(p^2) \times \mathcal{O}(p^2)$ and $B$ is $\mathcal{O}(p^2) \times \mathcal{O}(p)$. The almost-banded matrix $L$ may be written as the sum of an $\mathcal{O}(p)$-banded matrix $A$ and a rank-$\mathcal{O}(p)$ correction, $L = A + UCV^T$, where $U$ and $V$ are $\mathcal{O}(p^2) \times \mathcal{O}(p)$ and $C$ is $\mathcal{O}(p) \times \mathcal{O}(p)$. Using the Woodbury formula, the solution to $LX=B$ becomes
\[
X = L^{-1} B = \left(A + UCV^T\right)^{-1} B = \left( I \;-\; A^{-1} U \left( C^{-1} + V^T A^{-1} U \right)^{-1} V^T \right) A^{-1} B.
\]
The banded matrix $A$ can be inverted in $\mathcal{O}(p^3)$ operations and its inverse applied to $\mathcal{O}(p)$ righthand sides in $\mathcal{O}(p^4)$ operations. The matrix $C^{-1} + V^T A^{-1} U$ is $\mathcal{O}(p) \times \mathcal{O}(p)$ and so its inverse can be applied to $\mathcal{O}(p)$ righthand sides in $\mathcal{O}(p^3)$ operations. Therefore, the solution operator $S_{\mathcal{E}_i}$ on an element can be constructed in $\mathcal{O}(p^4)$ operations. The Dirichlet-to-Neumann operator $\dtn_{\mathcal{E}_i}$ can be computed as a matrix product in $\mathcal{O}(p^4)$ operations. As these operators are computed once for each element, the overall cost of the local build stage scales as
\[
\frac{p^4}{h^2} \approx N p^2.
\]

The cost of the global build stage and solve stage depends on the merge scheme defined by the indices $\mathcal{I}$. If the mesh $\mathcal{E} = \{\mathcal{E}_i\}_{i=1}^{n_\text{elem}}$ is approximately tensor-product, the merge indices $\mathcal{I}$ can be defined so that the hierarchy is approximately a binary tree (i.e., a binary tree with $\mathcal{O}(1)$ additional merges). If the mesh is unstructured, a hierarchical partitioning of the mesh may be computed by conversion to a graph partitioning problem~\cite{Karypis_98_01}. The partitioning should be as balanced as possible, so that the indices $\mathcal{I}$ define a balanced tree. If the user specifies merge indices that correspond to an unbalanced tree, then the tree may be automatically rebalanced. We assume that the merge indices $\mathcal{I}$ have been given so that the hierarchy in the build and solve stages approximately forms a binary tree with $\mathcal{O}(\log n_\text{elem})$ levels.

Let level $\ell = 0$ denote the bottom level of the hierarchy, where no elements have been merged. For a merge between $\mathcal{E}_i$ and $\mathcal{E}_j$ on level $\ell$ of the build stage, the solution operator $S_{\Gamma_{ij}}$ is computed by solving the linear system \cref{eq:build_sol}. The agglomerates $\mathcal{E}_i$ and $\mathcal{E}_j$ each contain $\mathcal{O}(2^\ell)$ mesh elements, with the interface between them, $\Gamma_{ij}$, containing $\mathcal{O}(2^{\ell/2})$ boundaries. Hence, the linear system in \cref{eq:build_sol} is $\mathcal{O}(2^{\ell/2} p) \times \mathcal{O}(2^{\ell/2} p)$ and can be solved in $\mathcal{O}((2^{\ell/2} p)^3)$ operations. As level $\ell$ has $\mathcal{O}(2^{-\ell} n_\text{elem})$ elements, the cost of processing all merges on level $\ell$ scales as
\[
\left(2^{-\ell} n_\text{elem}\right) \cdot \left(2^{\ell/2} p\right)^3 \;=\; n_\text{elem} 2^{\ell/2} p^3.
\]
The total cost for the global build stage then scales as
\[
p^3 n_\text{elem} \!\!\!\!\!\! \sum_{\ell=0}^{\mathcal{O}(\log n_\text{elem})}\!\!\!\! 2^{\ell/2}
\;\approx\; p^3 (n_\text{elem})^{3/2}
\;\approx\; \frac{p^3}{h^3}
\;\approx\; N^{3/2}
\]
as $N \to \infty$.

At level $\ell > 0$ of the solve stage, the unknown interface data is computed via a matrix-vector multiply with an $\mathcal{O}(2^{\ell/2}p) \times \mathcal{O}(2^{\ell/2}p)$ matrix. As level $\ell$ has $\mathcal{O}(2^{-\ell} n_\text{elem})$ elements, the cost of computing the solution on all interfaces scales as
\[
\left(2^{-\ell} n_\text{elem}\right) \cdot \left(2^{\ell/2}p\right)^2 \;=\; p^2 n_\text{elem}.
\]
The total cost for all levels $\ell > 0$ is then
\[
\sum_{\ell=1}^{\mathcal{O}(\log n_\text{elem})} p^2 n_\text{elem} \;\approx\; p^2 n_\text{elem} \log n_\text{elem}.
\]
At the bottom level, $\ell = 0$, the solution is computed on each element through matrix-vector multiplication with local solution operators of size $(p+1)^2 \times (4(p+1)+1)$, which requires $\mathcal{O}(p^3)$ operations. Therefore, the total cost for the solve stage scales as
\[
p^2 n_\text{elem} \log n_\text{elem} + p^3 n_\text{elem} \;\approx\; \frac{p^2}{h^2} \log \frac{1}{h^2} + \frac{p^3}{h^2} \;\approx\; N \log \tfrac{1}{h^2} + N p.
\]
The overall computational complexity of the method is therefore
\[
\underbrace{\vphantom{\tfrac12}Np^2 \,+\, \vphantom{\tfrac12}N^{3/2}}_{\text{build stage}} \;\;+\;\; \underbrace{N \log \tfrac{1}{h^2} + N p}_{\text{solve stage}} \;\approx\; Np^2 + N^{3/2}.
\]

As the method stores dense solution operators and Dirichlet-to-Neumann operators on every level of the hierarchy, the total storage cost is analogous to the computational cost of the solve stage. The amount of storage required by the method scales as
\[
N \log \tfrac{1}{h^2} + N p.
\]
The storage cost can become prohibitive when $p$ is large, as the local solution operators on each element require $\mathcal{O}(p^3 n_\text{elem})$ storage. However, these operators need not be constructed and stored. In the build stage, local Dirichlet-to-Neumann operators can be constructed directly by locally solving the PDE, evaluating the outward flux, and then discarding the solution. In the solve stage, the solution on the interior of each element can be computed by locally solving the PDE on the fly. This reduces the storage cost to $\mathcal{O}(N \log \tfrac{1}{h^2})$ while increasing the computational cost of the solve stage to $\mathcal{O}(N \log \tfrac{1}{h^2} + N p^2)$ operations, but does not change the overall computational complexity of the method.

\section{Software}\label{sec:software}
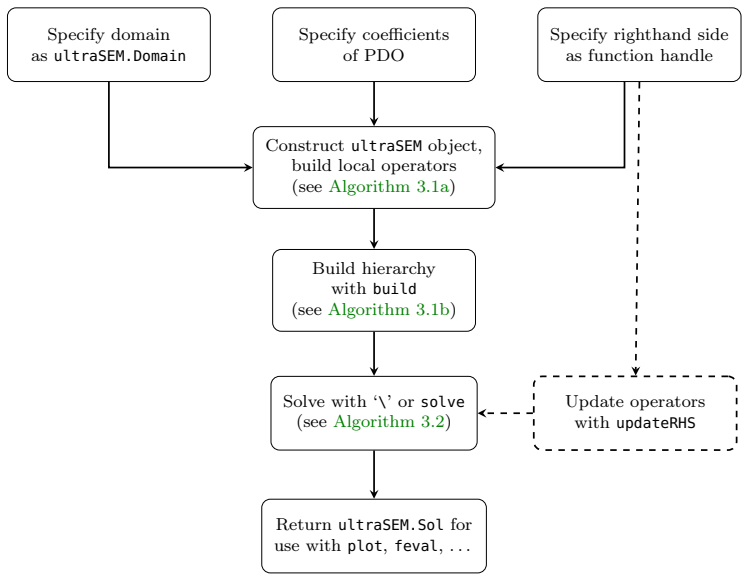
\begin{figure}[htb]
	\centering
	\vspace{0.2cm}
	\resizebox{0.75\textwidth}{!}{%
	\begin{tikzpicture}[node distance=2cm, every node/.style={minimum width=3.3cm, minimum height=1.2cm, draw=black, inner sep=0.2cm, align=center, font=\footnotesize}]
		\tikzstyle{arrow} = [thick,->,>=stealth]
		\node (domain) [rectangle, rounded corners] {Specify domain \\ as \mytt{ultraSEM.Domain}};
		\node (pdo) [rectangle, rounded corners, right=1cm of domain] {Specify coefficients \\ of PDO};
		\node (rhs) [rectangle, rounded corners, right=1cm of pdo] {Specify righthand side \\ as function handle};
		\node (construct) [rectangle, rounded corners, below of=pdo] {Construct \mytt{ultraSEM} object, \\ build local operators \\ (see \cref{alg:init})};
		\draw [arrow] (domain) |- (construct);
		\draw [arrow] (pdo) -- (construct);
		\draw [arrow] ([xshift=-8mm]rhs) |- (construct);
		\node (build) [rectangle, rounded corners, below of=construct] {Build hierarchy \\ with \mytt{build} \\ (see \cref{alg:build})};
		\draw [arrow] (construct) -- (build);
		\node (solve) [rectangle, rounded corners, below of=build] {Solve with `\mytt{\char`\\}' or \mytt{solve} \\ (see \cref{alg:solve})};
		\draw [arrow] (build) -- (solve);
		\node (sol) [rectangle, rounded corners, below of=solve] {Return \mytt{ultraSEM.Sol} for \\ use with \mytt{plot}, \mytt{feval}, \textellipsis};
		\draw [arrow] (solve) -- (sol);
		\node (updateRHS) [rectangle, rounded corners, right=0.9cm of solve, thick, dashed] {Update operators \\ with \mytt{updateRHS}};
		\draw [arrow, dashed] (updateRHS.west) -- (solve.east);
		\draw [arrow, dashed] (rhs.south) -- (updateRHS);
	\end{tikzpicture}%
	}
	\vspace{0.2cm}
	\caption{A diagram of the code workflow in \ultraSEM. The code is designed to mirror the steps of the hierarchical Poincar\'{e}--Steklov scheme.}
	\label{fig:code_flow}
\end{figure}

\begin{figure}[htb]
\centering
    \begin{minipage}{0.33\textwidth}
    \centering
    \includegraphics[width=0.67\textwidth]{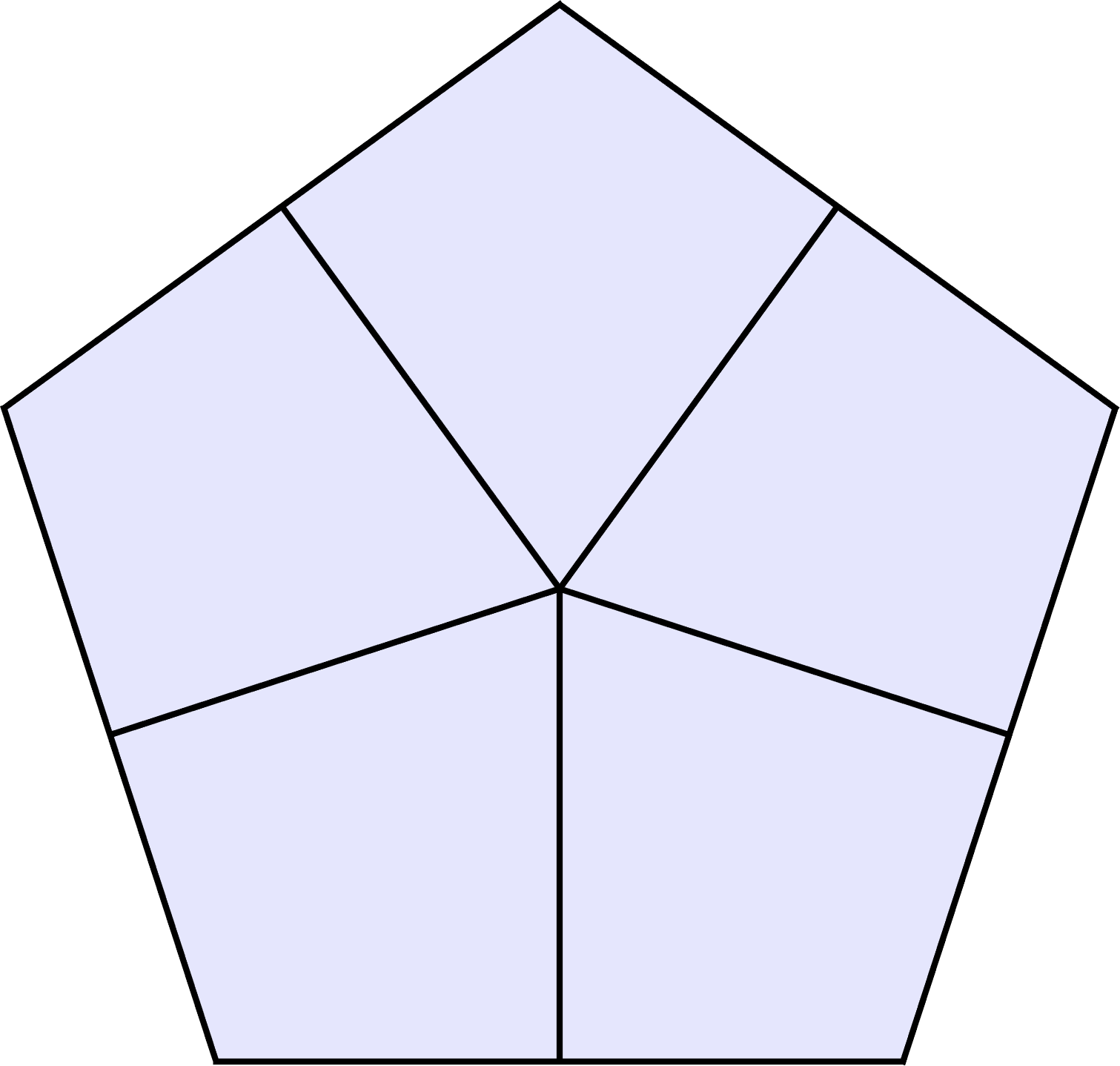} \\
    \begin{lstlisting}[style=Matlab-editor, basicstyle=\footnotesize\mlttfamily]
   dom = ultraSEM.polygon(5);
   plot(dom)
    \end{lstlisting}
    \vspace{1cm}
    \end{minipage}%
    \hspace{-0.3cm}%
    \begin{minipage}{0.33\textwidth}
    \centering
    \includegraphics[width=0.67\textwidth]{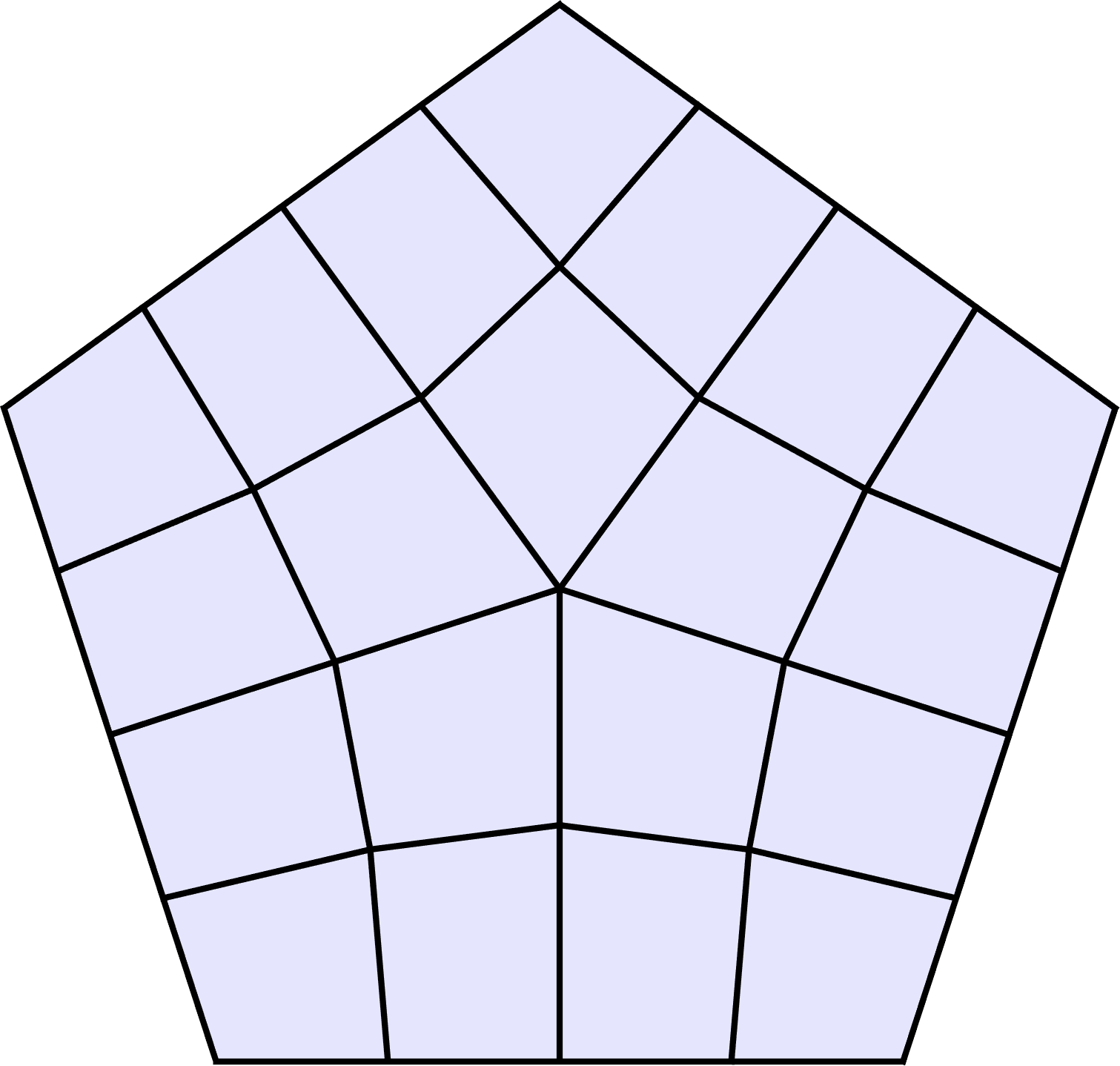} \\
    \begin{lstlisting}[style=Matlab-editor, basicstyle=\footnotesize\mlttfamily]
      dom = refine(dom);
      plot(dom)
    \end{lstlisting}
    \vspace{1.05cm}
    \end{minipage}%
    \hspace{-0.1cm}%
    \begin{minipage}{0.33\textwidth}
    \centering
    \vspace{-0.15em}
    \includegraphics[width=0.67\textwidth]{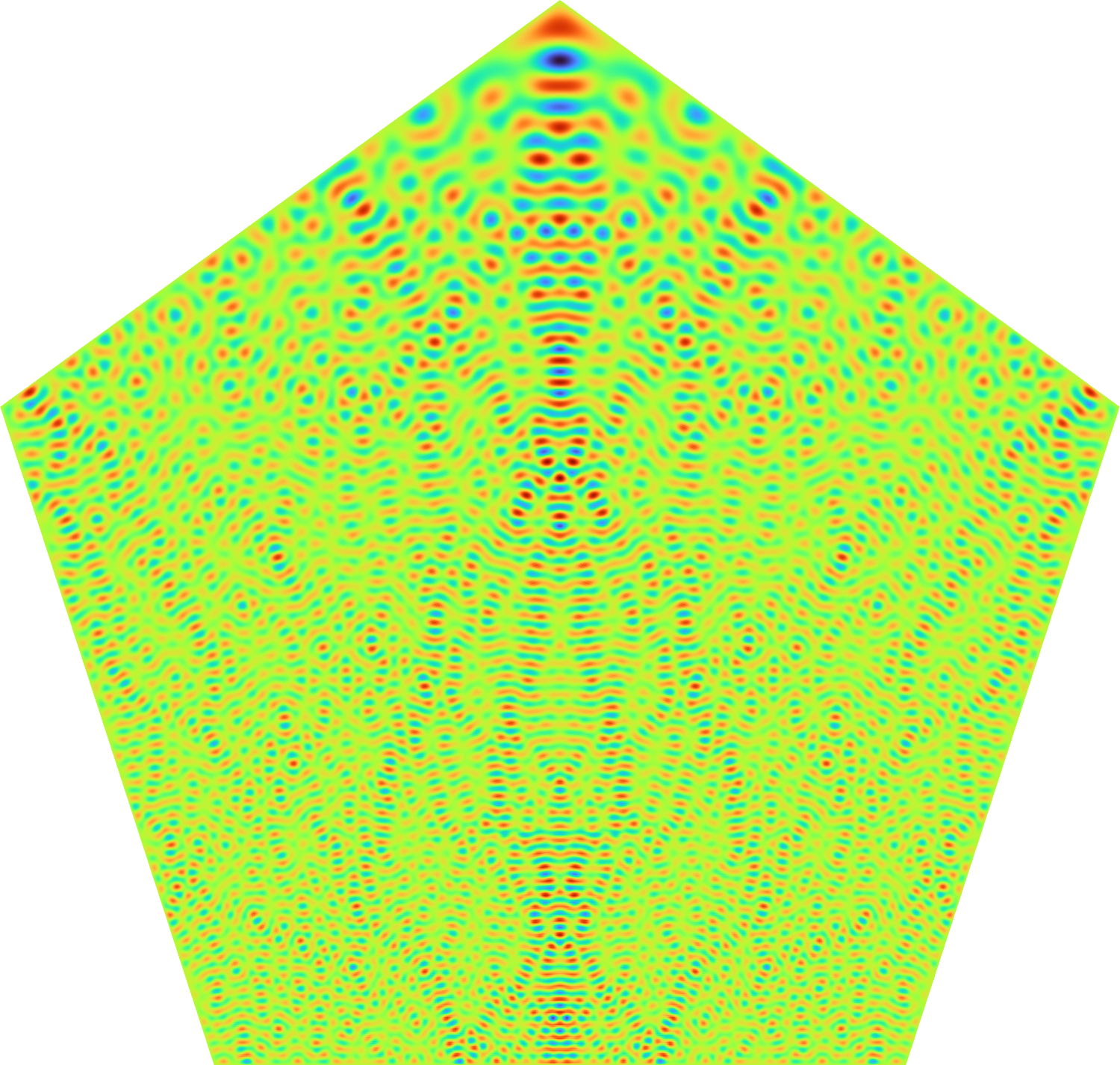}
    \begin{lstlisting}[style=Matlab-editor, basicstyle=\footnotesize\mlttfamily, linewidth=4.6cm]
 p = 60; rhs = -1; bc = 0;
 pdo = {{1,0,1}, {0,0}, ...
        @(x,y) 50000*(1-y)};
 S = ultraSEM(dom, pdo, rhs, p);
 u = S \ bc; plot(u)
    \end{lstlisting}
    \end{minipage}
    \vspace{-0.8em}
\caption{A simple example of the syntax in \ultraSEM. A pentagonal domain (with side length 1.2) is meshed into five quadrilaterals (left) and uniformly refined (center). The gravity Helmholtz equation $\nabla^2 u + 50000(1-y)u = -1$ with zero Dirichlet boundary conditions is then solved on the mesh using polynomials of degree 60 on each element (right).}
\label{fig:code_ex}
\end{figure}

We have implemented the ultraspherical SEM in an open-source software package, \ultraSEM, written in MATLAB without parallelization~\cite{GithubRepo}. An outline of the workflow is depicted in \cref{fig:code_flow}, and a simple example is shown in \cref{fig:code_ex}.

The user constructs each element as an \mytt{ultraSEM.Domain}, which encodes the coordinate transformations and merge indices local to each element. Convenient functions for constructing rectangles, quadrilaterals, triangles, and polygons are available via the commands \mytt{ultraSEM.rectangle}, \mytt{ultraSEM.quad}, \mytt{ultraSEM.triangle}, and \mytt{ultraSEM.polygon}, respectively (see \cref{fig:code_ex} (left)), which automatically encode the suitable transformations and merge indices. Elements can be combined to form larger domains by merging them with the `\mytt{\&}' operator; the merge indices $\mathcal{I}$ will then correspond to the order induced by the sequence of `\mytt{\&}' operations. More general meshes can be constructed using the \mytt{refine(dom)} method (see \cref{fig:code_ex} (center)), which performs uniform $h$-refinement on a given domain \mytt{dom}, or the \mytt{refinePoint(dom, [x,y])} method, which performs adaptive $h$-refinement on \mytt{dom} around the point $(x,y)$ (see \cref{fig:corner_point_refinement}).

A PDO is specified by its coefficients for each derivative, in the form \mytt{\{\{uxx, uxy, uyy\}, \{ux, uy\}, b\}}, where each term \mytt{uxx}, \mytt{uxy}, \textellipsis\ can be a scalar (constant coefficient) or function handle (variable coefficient). The domain and PDO are then passed---along with a righthand side and polynomial order---to construct an \mytt{ultraSEM} object (see \cref{fig:code_ex} (right)). The \mytt{ultraSEM} constructor initializes the local operators on each element (see \cref{alg:init}), which are represented as \mytt{ultraSEM.Leaf} objects in the hierarchy. The hierarchy of merged operators may then be built in an upward pass via the \mytt{build} command (see \cref{alg:build}), which creates a tree of \mytt{ultraSEM.Parent} objects (if \mytt{build} is not explicitly called, the build stage is automatically performed when the user requests a solve to be executed). The solve stage is invoked via the \mytt{solve} command (or equivalently, the `\mytt{\char`\\}' operator), which computes the solution by applying the hierarchy of operators in a downward pass (see \cref{alg:solve}). The solution is returned as an \mytt{ultraSEM.Sol} object, which overloads a host of functions for plotting (e.g., \mytt{plot}, \mytt{contour}) and evaluation (e.g., \mytt{feval}, \mytt{norm}).

An \mytt{ultraSEM} object that has been built can be repeatedly applied to new boundary conditions by invoking \mytt{solve} multiple times. The object can also be cheaply updated to solve with a new righthand side by calling \mytt{updateRHS}, which alters the last column of each operator in the hierarchy to correspond to a new particular solution.

\section{Numerical results}\label{sec:results}

\subsection{Computational complexity}
\begin{figure}[htb]
  \centering
  \begin{overpic}[width=0.48\textwidth]{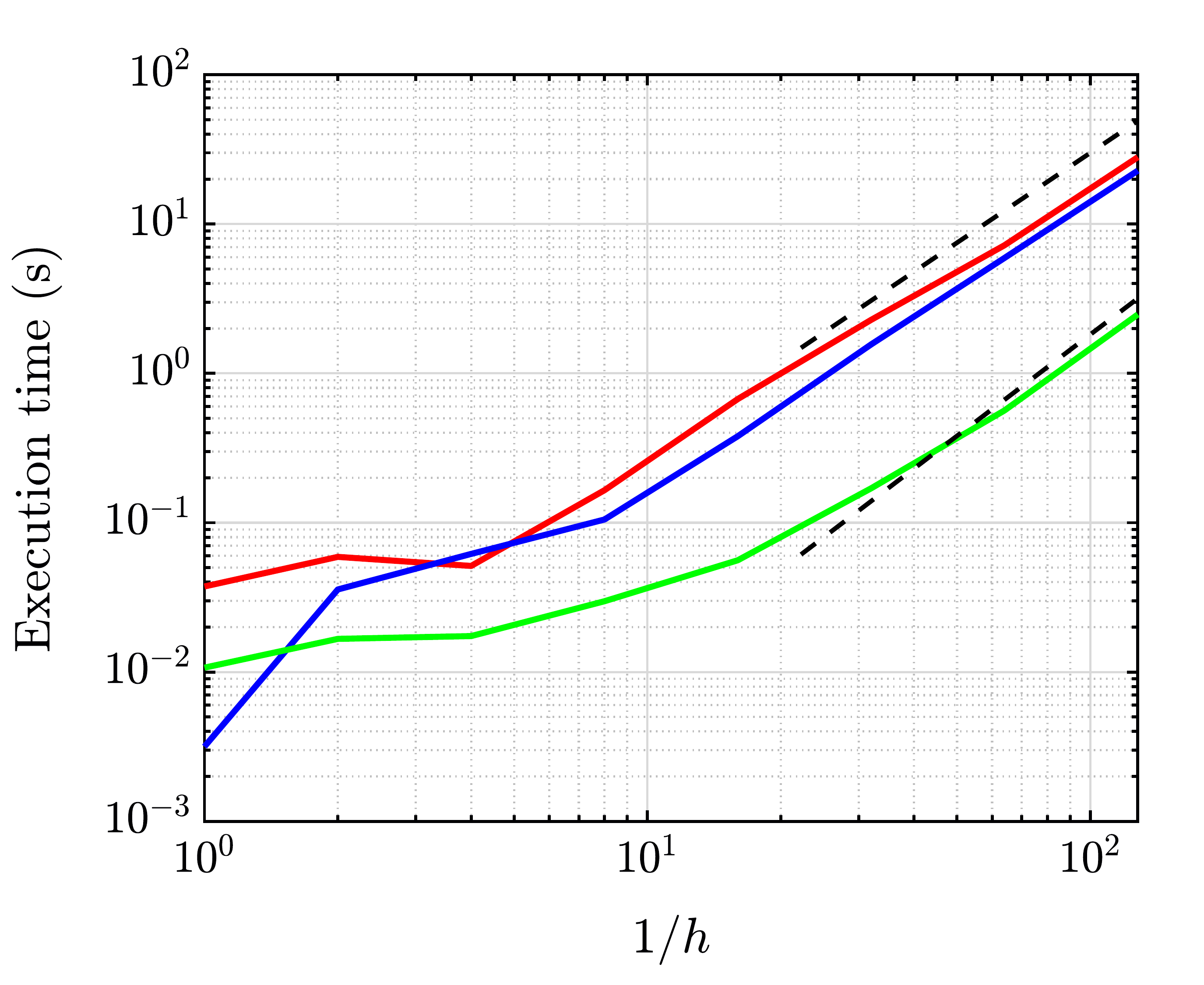}%
    \put(73,60.5) {\rotatebox{33}{\scalebox{0.65}{$\mathcal{O}(1/h^2)$}}}
    \put(66,38.8) {\rotatebox{37}{\scalebox{0.65}{$\mathcal{O}((1/h^2) \log(1/h^2))$}}}
  \end{overpic}%
  ~~~~%
  \begin{overpic}[width=0.48\textwidth]{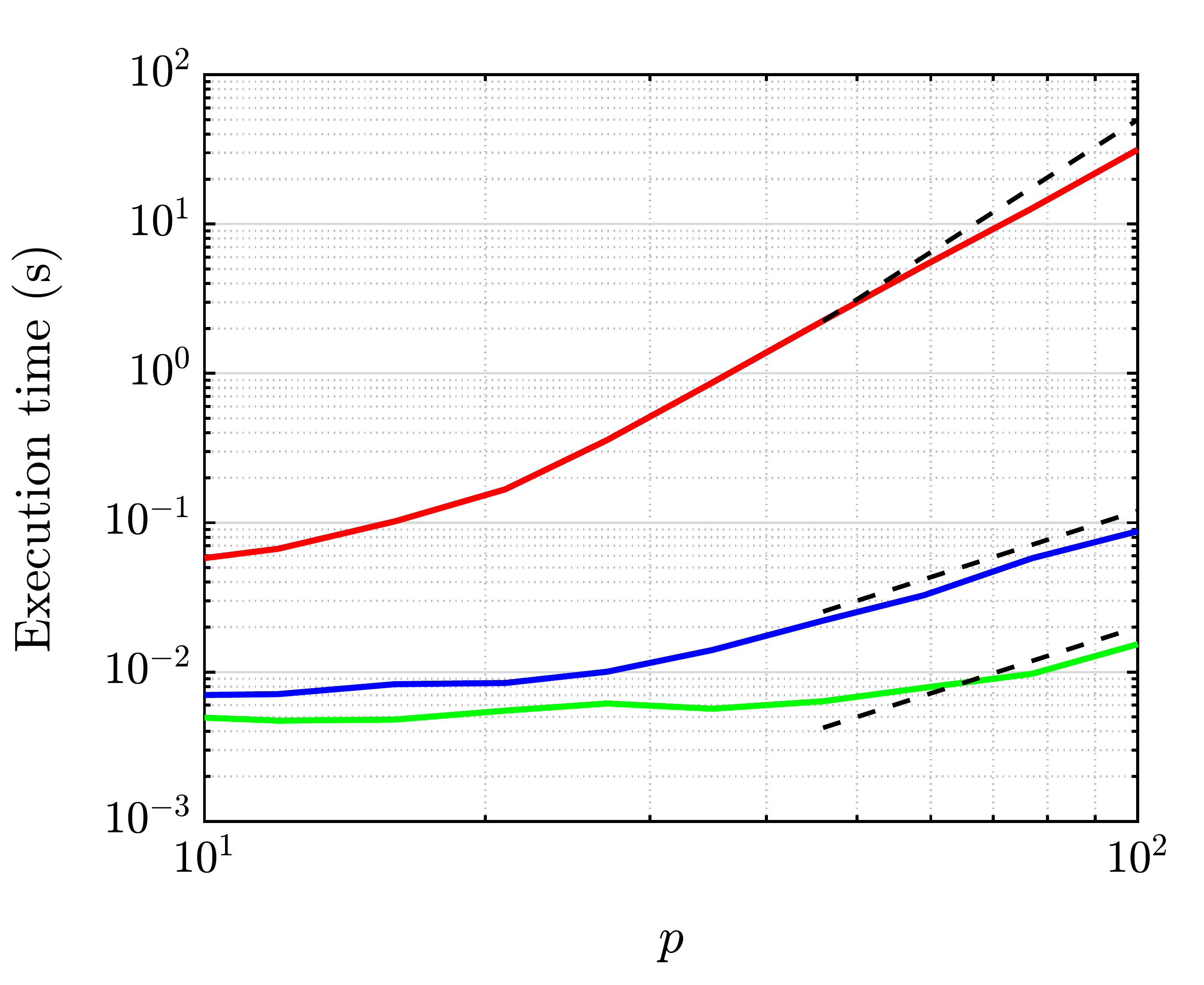}%
    \put(76,63) {\rotatebox{31}{\scalebox{0.65}{$\mathcal{O}(p^4)$}}}
    \put(76,35.5) {\rotatebox{17}{\scalebox{0.65}{$\mathcal{O}(p^2)$}}}
    \put(77,26) {\rotatebox{17}{\scalebox{0.65}{$\mathcal{O}(p^2)$}}}
  \end{overpic}%
  \caption{Execution time (in seconds) for \ultraSEM over a range of mesh sizes (left) and polynomial orders (right). Timings are depicted for the local build stage (red), global build stage (blue), and solve stage (green), when solving the PDE $\nabla^2 u + \sin(xy) u = f$ on the domain $\Omega = [0,1]^2$ with spatially varying righthand side and spatially varying Dirichlet boundary conditions. On the left, we successively refine a Cartesian mesh while keeping the polynomial order fixed at $p=4$. On the right, we use a $4 \times 4$ Cartesian mesh while successively increasing the polynomial order.}
  \label{fig:complexity}
\end{figure}
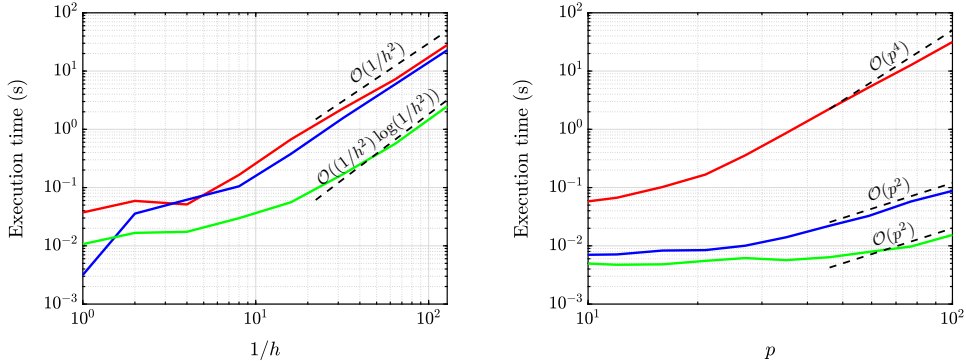

To illustrate the computational complexity of \ultraSEM, we measure the execution times\footnote{All numerical experiments were performed in MATLAB R2020a on a 40-core Intel Xeon E5-2630 workstation with 128GB of RAM and no explicit parallelization.} of the build and solve stages of the method under uniform $h$- and $p$-refinement. \cref{fig:complexity} shows the recorded timings for solving the variable coefficient PDE $\nabla^2 u + \sin(xy) u = f$ on the domain $\Omega = [0,1]^2$ with a spatially varying righthand and spatially varying Dirichlet boundary conditions.

In \cref{fig:complexity} (left), the polynomial order is fixed at $p=4$ and a Cartesian mesh with $\mathcal{O}(1/h^2)$ elements is successively refined. Both the local and global build stages exhibit $\mathcal{O}(1/h^2)$ scaling as $h \to 0$, while the solve stage scales as $\mathcal{O}(1/h^2 \log(1/h^2))$. The timings for the build stage do not exhibit the expected $\mathcal{O}(1/h^3)$ scaling. This is likely due to the fact that the build stage relies on dense linear algebra routines that have been heavily optimized for the relatively small $\mathcal{O}(1/h) \times \mathcal{O}(1/h)$ matrices tested here. The storage used by the finest mesh in \cref{fig:complexity} (left) is approximately 2GB.

In \cref{fig:complexity} (right), the Cartesian mesh is fixed to have $4 \times 4$ elements and the polynomial order $p$ is successively increased. The cost of the local build stage dominates, exhibiting close to the expected $\mathcal{O}(p^4)$ scaling as $p \to \infty$. The global build and solve stages perform better than expected, both exhibiting $\mathcal{O}(p^2)$ scaling. Again, this can likely be attributed to the performance of dense linear algebra routines in the regime of $p$ tested. The storage used by the finest mesh in \cref{fig:complexity} (right) is approximately 25GB.

\subsection{Convergence and $hp$-adaptivity}
\begin{figure}[htb]
  \centering
  \begin{overpic}[width=0.48\textwidth]{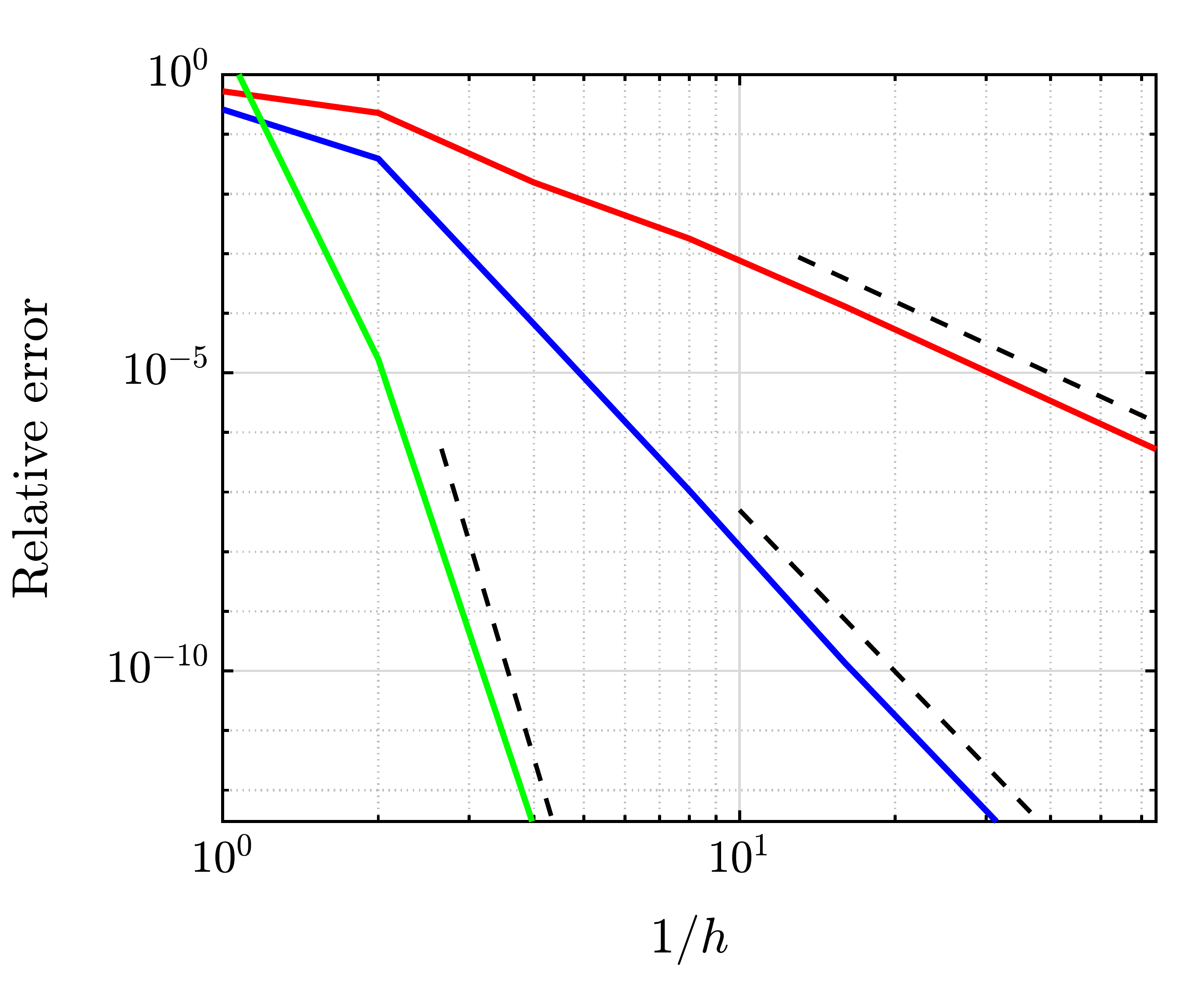}%
    \put(73,58.9) {\rotatebox{-24}{\scalebox{0.65}{$p=5$,\;\;$\mathcal{O}(h^4)$}}}
	\put(67.5,35) {\rotatebox{-46}{\scalebox{0.65}{$p=10$,\;\;$\mathcal{O}(h^9)$}}}
	\put(39.1,39.4) {\rotatebox{-72.5}{\scalebox{0.65}{$p=30$,\;\;$\mathcal{O}(h^{29})$}}}
  \end{overpic}%
  ~~~~%
  \begin{overpic}[width=0.48\textwidth]{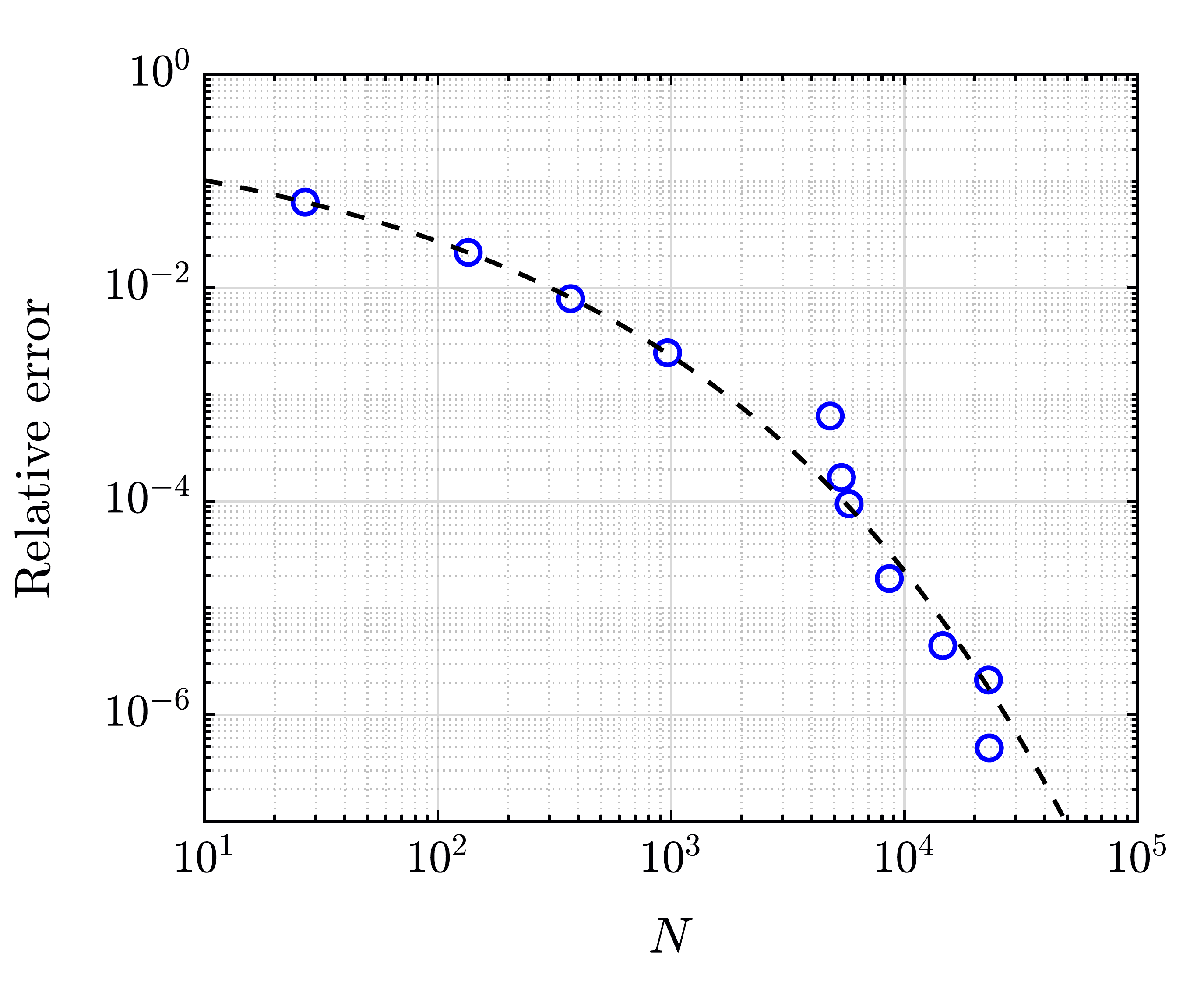}%
    \put(43,65) {\rotatebox{-34}{\scalebox{0.75}{$\mathcal{O}\!\left(e^{-0.8N^{0.27}}\right)$}}}
  \end{overpic}
  \caption{Convergence of \ultraSEM with respect to $h$ and $p$. (Left) Relative error in the $L^2$ norm when solving \cref{eq:hp_convergence1} with $\omega = p$ under uniform $h$-refinement, for $p=5$ (red), $p=10$ (blue), and $p=30$ (green). In each case, $\mathcal{O}(h^{p-1})$ convergence is observed. (Right) An \emph{a priori} $hp$-adaptivity strategy is applied to the L-shape problem in \cref{eq:Lshape}. The relative error decays super-algebraically in the total number of degrees of freedom $N$.}
  \label{fig:hp_convergence}
\end{figure}

We now investigate the convergence properties of \ultraSEM with respect to the mesh size $h$ and polynomial order $p$. As a test problem, we consider solving the Helmholtz equation,
\begin{equation}\label{eq:hp_convergence1}
\nabla^2 u + (\sqrt{2} \omega)^2 u = 0, \qquad u \in [-1,1]^2,
\end{equation}
with $\omega \in \mathbb{R}$ and Dirichlet boundary conditions given so that the exact solution is $u(x,y) = \cos(\omega x) \cos(\omega y)$. To measure convergence over a range of polynomial orders, we set $\omega = p$ so that the number of degrees of freedom per wavelength remains fixed independent of $p$. We then solve \cref{eq:hp_convergence1} under uniform $h$-refinement. \cref{fig:hp_convergence} shows the relative error in the $L^2$ norm as $h \to 0$ for polynomial orders $p=5$, $p=10$, and $p=30$. The convergence rate is observed to be $\mathcal{O}(h^{p-1})$. If error is measured in the $H^1$ or $H^2$ norm, where $H^k$ denotes the Sobolev space of functions whose weak derivatives up to order $k$ are in $L^2$, then the convergence rate is similarly $\mathcal{O}(h^{p-1})$. Since our method is sparse with respect to $p$, the exact rate of convergence is not so important, as a degree-$p$ discretization may easily be replaced by a degree-$(p+1)$ discretization with minimal increase in computational cost.

\begin{figure}[htb]
  \centering
  \includegraphics[width=0.1\textwidth]{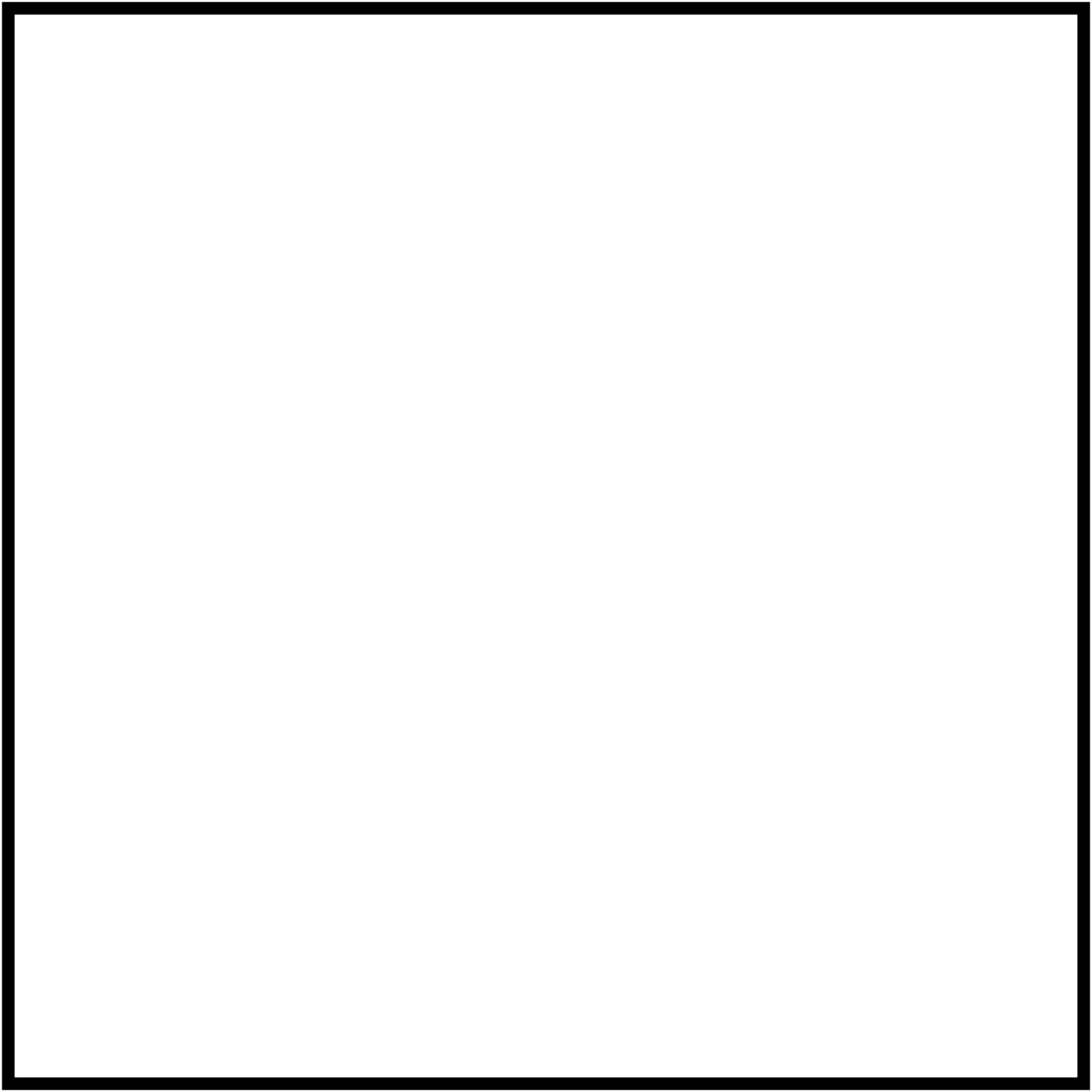}
  \hspace{0.4cm}
  \includegraphics[width=0.1\textwidth]{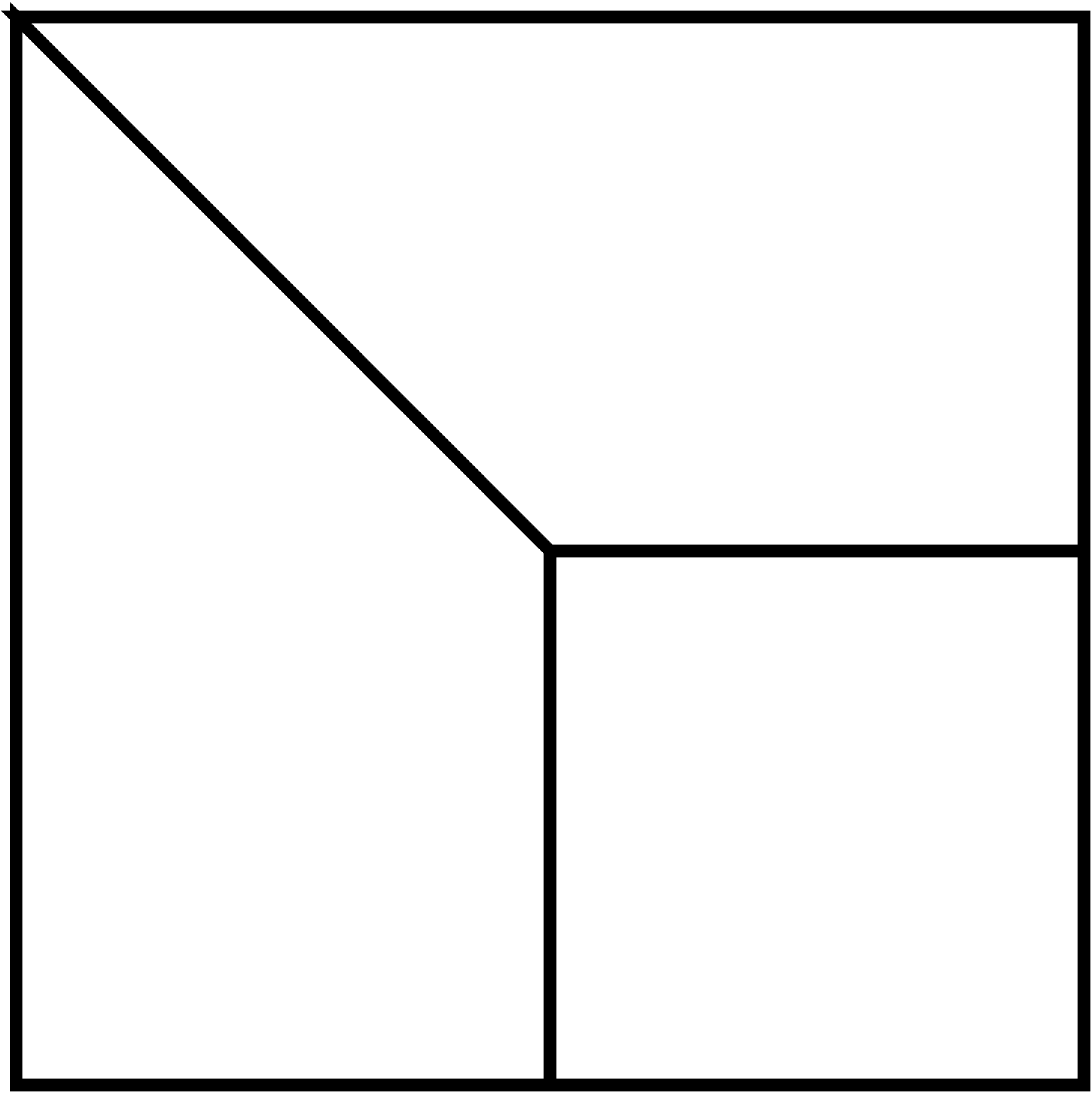}
  \hspace{0.4cm}
  \includegraphics[width=0.1\textwidth]{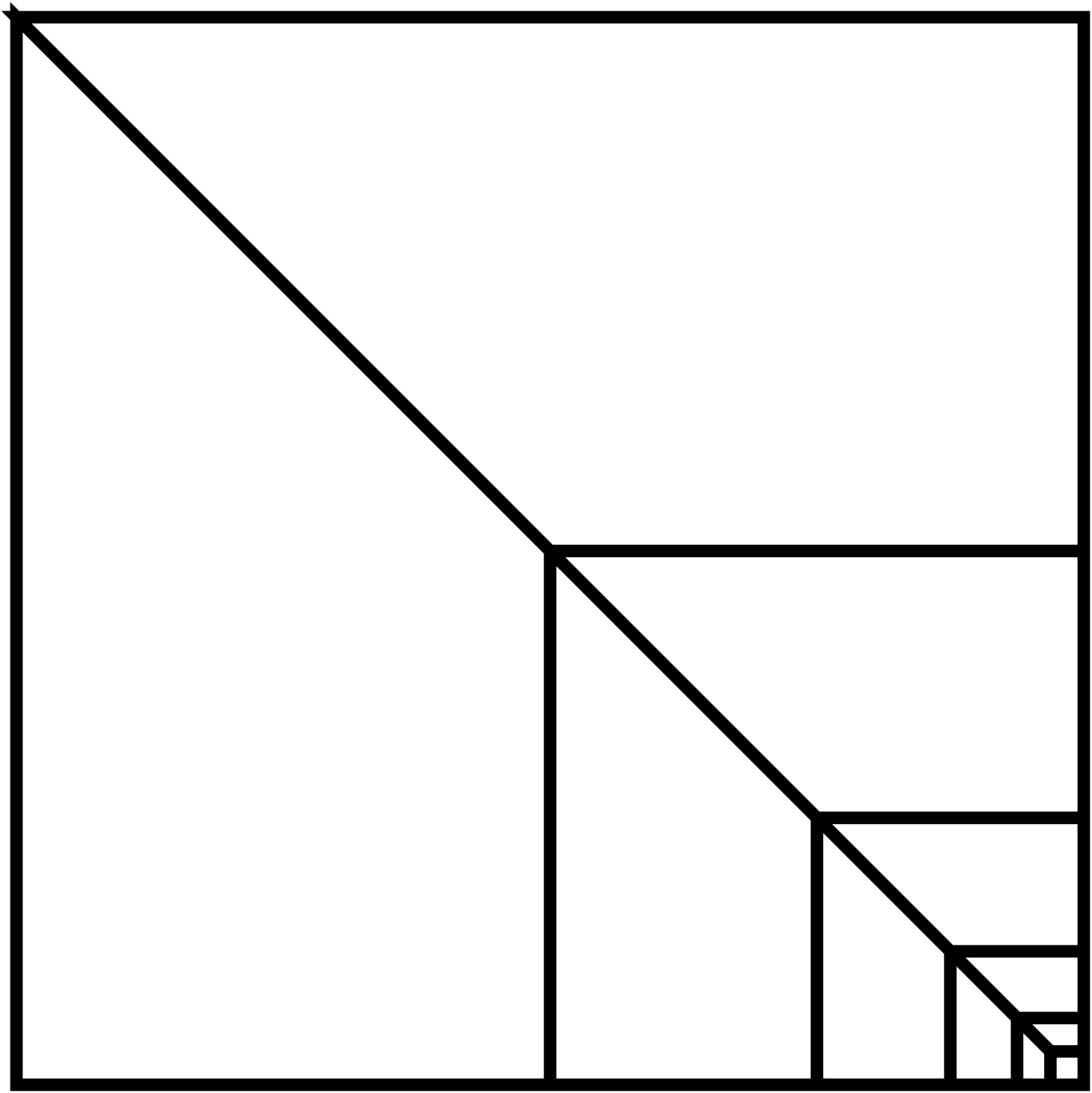}
  \hspace{1.5cm}
  \includegraphics[width=0.1\textwidth]{square}
  \hspace{0.4cm}
  \includegraphics[width=0.1\textwidth]{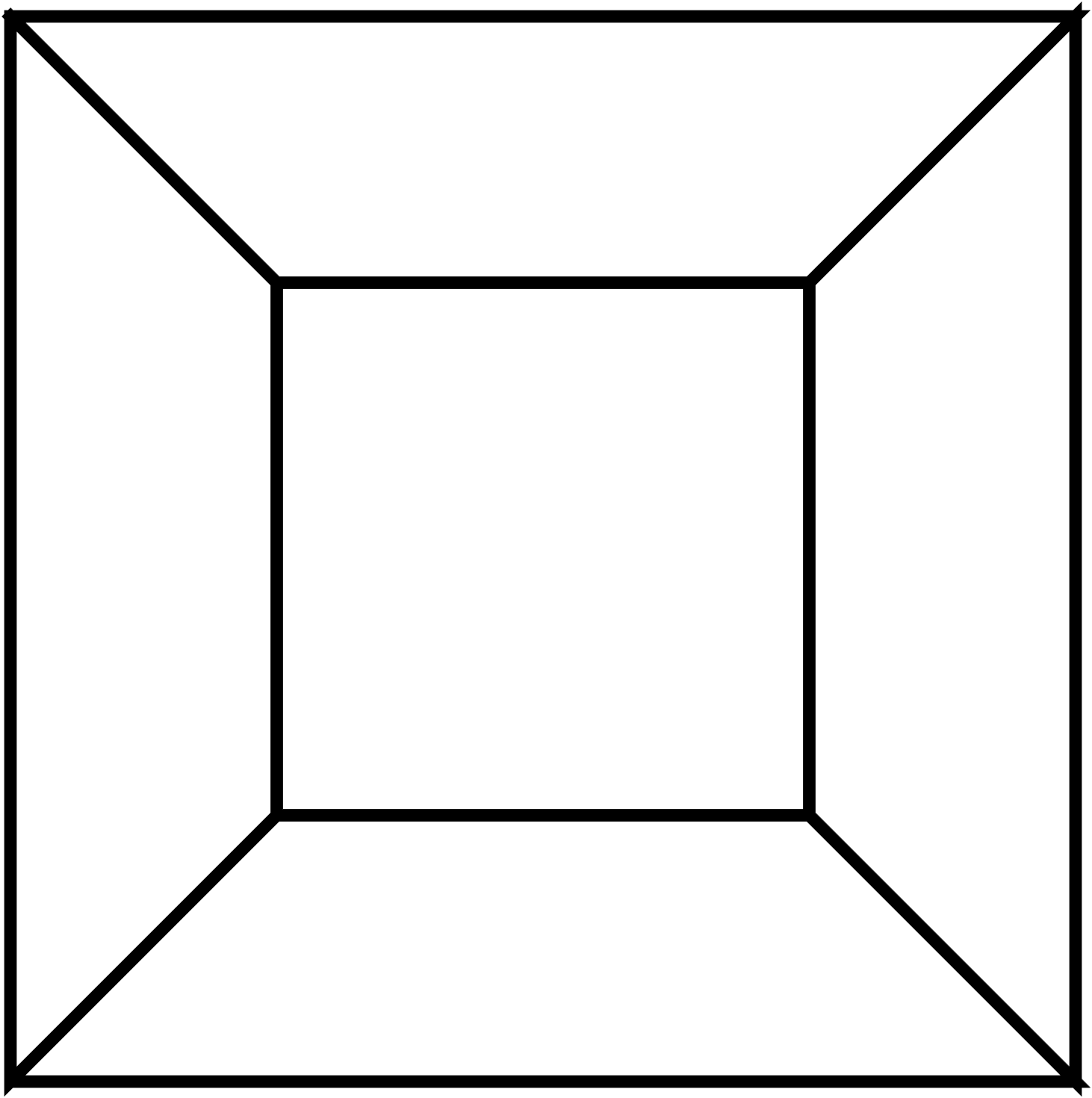}
  \hspace{0.4cm}
  \includegraphics[width=0.1\textwidth]{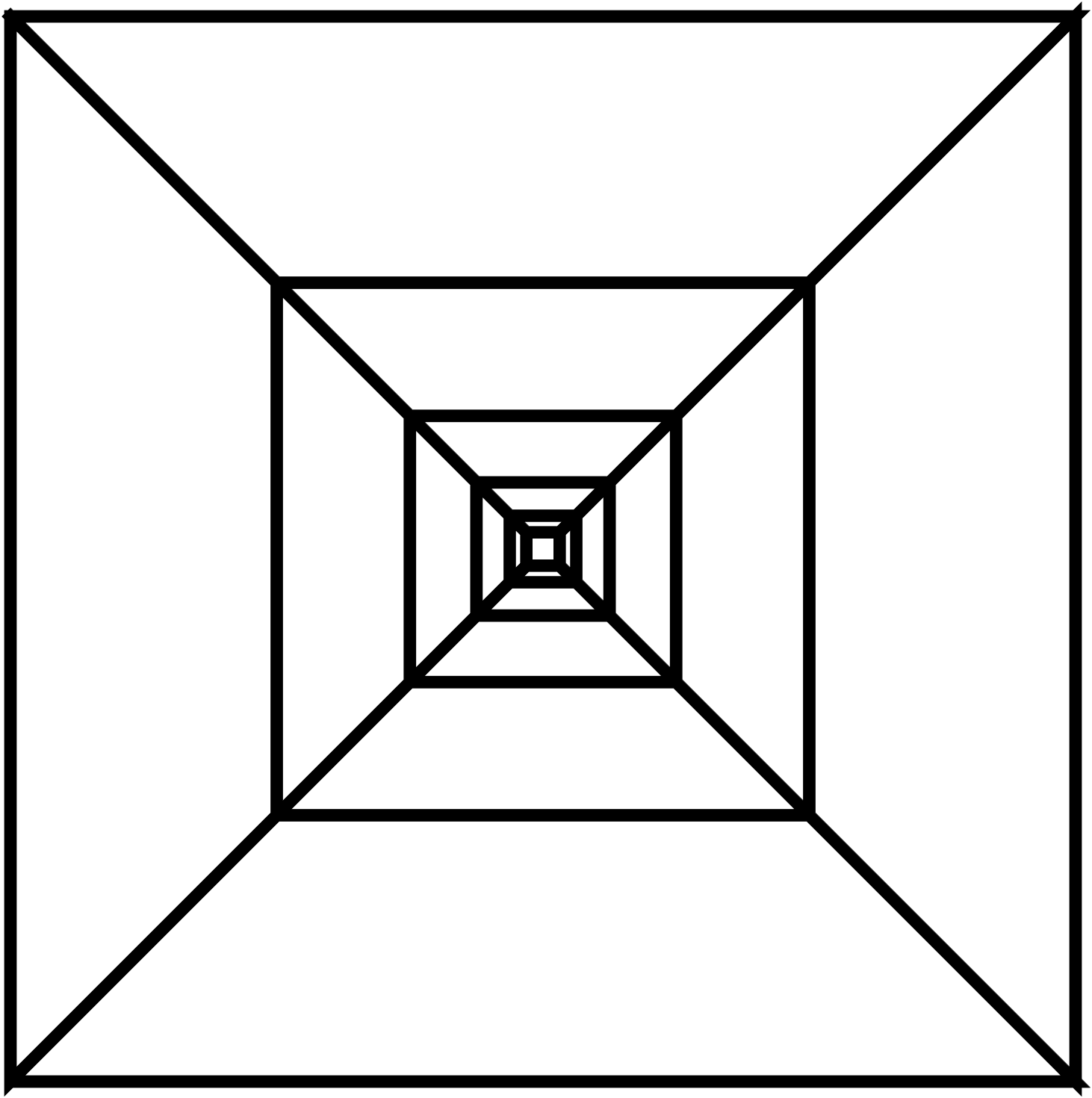}
  \vspace{0.3em}
  \caption{To avoid hanging nodes, \ultraSEM performs $h$-refinement in a conforming way. (Left)  A square is successively refined into a corner by subdividing into three children per refinement level. (Right) A square is successively refined around a point by subdividing into five children per refinement level. General quadrilaterals are refined in similar ways.}
  \label{fig:corner_point_refinement}
\end{figure}

In general, the mesh size $h$ and polynomial order $p$ need not be the same on each element. Adaptive $h$-refinement can be performed on each element locally; however, subdividing an element may give rise to meshes with hanging nodes (i.e., nodes of the mesh which occur in the middle of an element's face). While hanging nodes may be handled in the hierarchical Poincare--Steklov scheme through the use of interpolation operators~\cite{Geldermans_19_01}, we choose to avoid them here. To avoid hanging nodes, \ultraSEM performs $h$-refinement in a conforming way around specified corners or points, by subdividing a quadrilateral element into three or five children, respectively (see \cref{fig:corner_point_refinement}).

The ultraspherical spectral element method can naturally perform $p$-adaptivity by applying local interpolation and restriction operators to the elemental matrices involved in each merge operation. Since each unknown interface function is represented by a vector of Chebyshev coefficients, interpolation to and restriction from an interface can be performed simply by zero-padding or truncating the interface data. The polynomial order on an interface can be defined in a variety of ways. Popular choices include the minimum rule and maximum rule~\cite{Demkowicz_06_01}; we employ the minimum rule here, which sets the polynomial order on an interface to be the minimum of the polynomial orders on the adjacent elements.

We now consider the application of an $hp$-adaptivity strategy to the classical L-shape domain problem~\cite{Mitchell_13_01},
\begin{equation}\label{eq:Lshape}
\nabla^2 u = 0, \qquad u \in [-1,1]^2 \;\setminus\; [0,1] \times [-1,0],
\end{equation}
with Dirichlet boundary conditions given so that the exact solution is $u(r,\theta) = r^{2/3} \sin(2\theta/3)$, where $r = \sqrt{x^2 + y^2}$ and $\theta = \tan^{-1}(y/x)$. The reentrant corner of the domain induces a singularity in the solution so that $u \in H^{1+2/3}$ near the origin. Therefore, any strategy based on uniform $h$- or $p$-refinement is necessarily restricted to algebraic convergence\footnote{For this Laplace problem, alternative methods may provide higher accuracy per degree of freedom than element methods. For instance, root-exponential convergence in the supremum norm can be achieved by representing the solution as the real part of a rational function with poles exponentially clustered near each corner~\cite{Gopal_19_01}.}~\cite{Babuska_87_01}. That is, for a numerical solution $u_{hp}$ based on uniform refinement, the error can be bounded \emph{a priori} by
\[
\|u - u_{hp}\|_{L^2} \leq \|u - u_{hp}\|_{H^1} \leq C \left(\frac{h}{p}\right)^{2/3} \|u\|_{H^{1+2/3}},
\]
for some constant $C>0$. However, by employing a suitable $hp$-adaptivity strategy, super-algebraic convergence in the number of degrees of freedom $N$ can be achieved~\cite{Babuska_86_01}, i.e.,
\[
\|u - u_{hp}\|_{L^2} \leq \|u - u_{hp}\|_{H^1} \leq C_1 e^{-C_2 N^{1/3}},
\]
for some constants $C_1, C_2 > 0$. Here we employ an \emph{a priori} adaptivity strategy, where $h$-refinement is performed into the reentrant corner on elements adjacent to the origin and $p$-refinement is performed on all other elements~\cite{Ainsworth_98_01}. Given a desired relative error tolerance and an initial coarse $hp$-mesh, an automatic $hp$-adaptivity loop is run that successively refines or coarsens each element in $h$ or $p$ based on an \emph{a posteriori} error indicator~\cite{Mitchell_11_01, Ainsworth_97_01}. Here, we compute the element-wise error from the exact solution as a surrogate for a true \emph{a posteriori} error indicator. \cref{fig:hp_convergence} (right) shows the relative error in the $L^2$ norm versus the total number of degrees of freedom $N$ in the adaptive $hp$-mesh for a sequence of error tolerances. Super-algebraic convergence to the solution is observed as the number of degrees of freedom increases. A least-squares fit to the data gives an approximate convergence rate of $\mathcal{O}(e^{-0.8 N^{0.27}})$. To illustrate the range of $h$ and $p$ used on a given mesh, for a relative error tolerance of $10^{-6}$ the final mesh contains 15 levels of corner $h$-refinement and polynomial orders ranging from 3 to 13.

\begin{figure}[htb]
\centering
    \begin{minipage}{0.49\textwidth}
    \centering
    \begin{overpic}[width=0.8\textwidth]{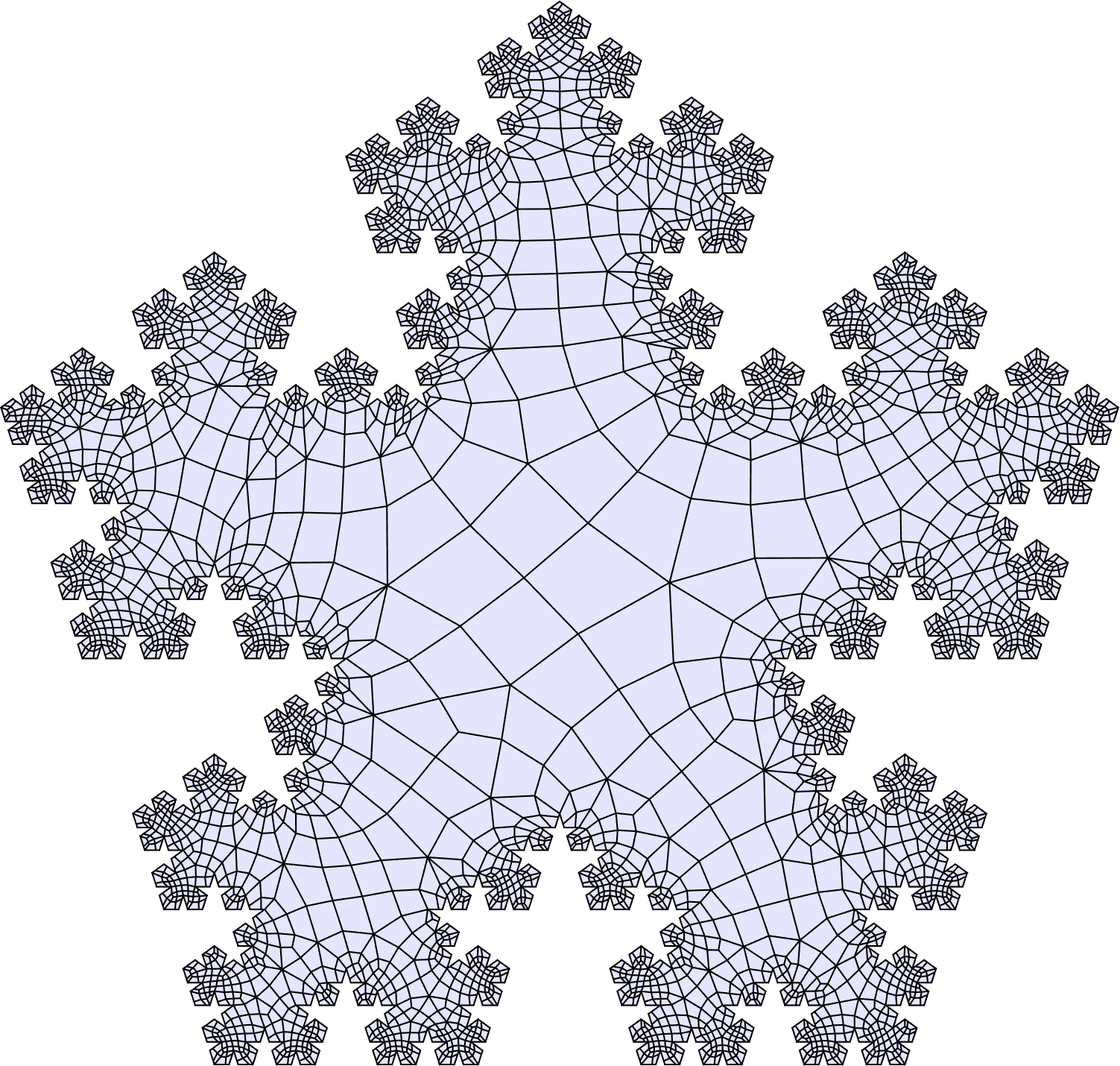}
        \put(89.5,8) {\rotatebox{-16.5}{\footnotesize 1.2}}
        \put(83.5,-1) {\begin{tikzpicture}\draw[{Bar[width=1mm]}-{Bar[width=1mm]}] (0,0) -- (0.35,1.15);\end{tikzpicture}}
    \end{overpic} \\
    \end{minipage}%
    \begin{minipage}{0.49\textwidth}
    \centering
    \includegraphics[width=0.8\textwidth]{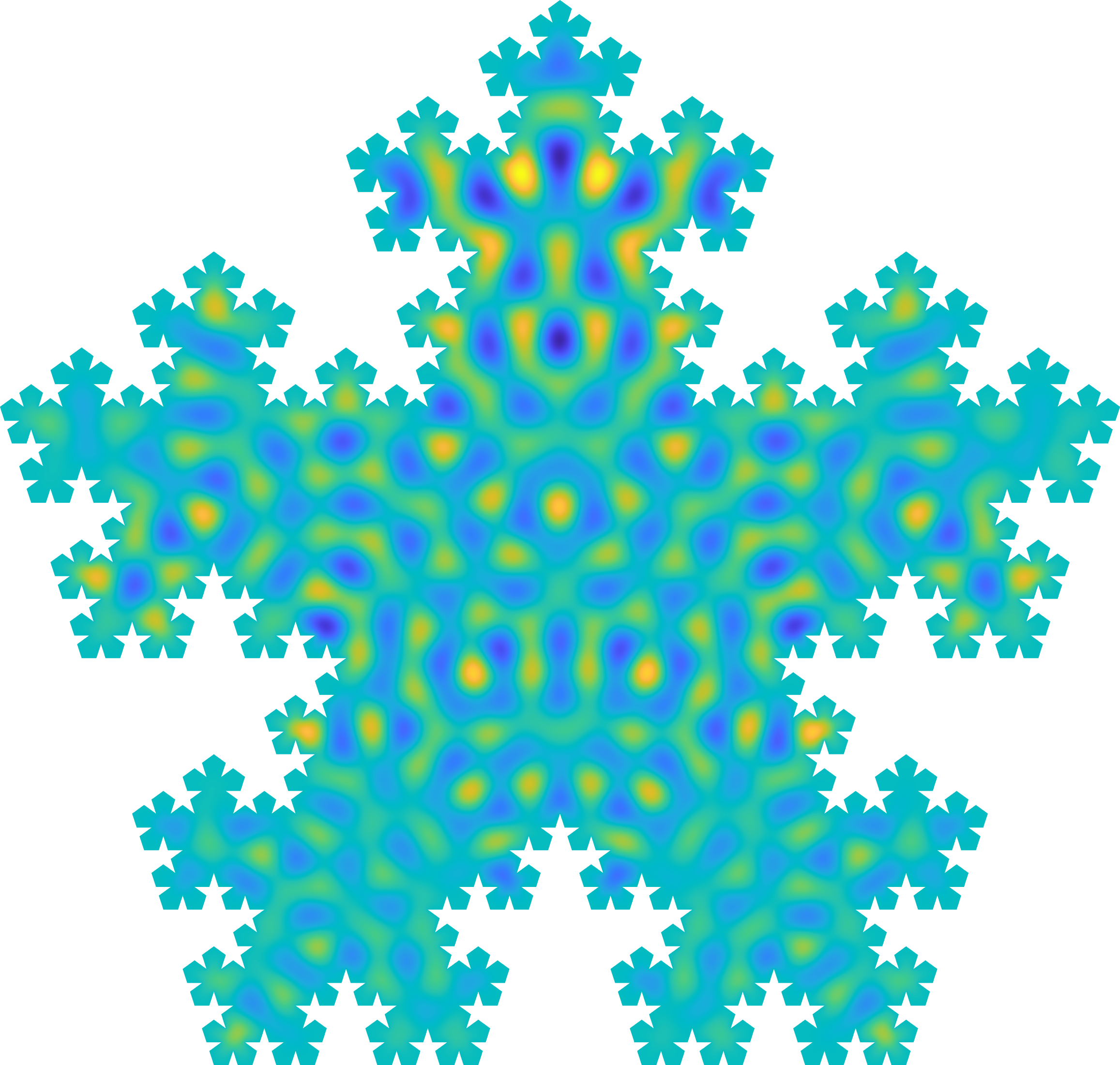} \\
    \end{minipage}%
\caption{(Left) A snowflake-shaped domain created by a Penrose tiling is adaptively meshed with 4,658 elements, with $p$ varying from $p=7$ on small elements to $p=40$ on large elements. (Right) The gravity Helmholtz equation \cref{eq:gravity_helmholtz} is solved on this domain. The solution is represented by $N = \text{333,627}$ degrees of freedom.}
\label{fig:snowflake}
\end{figure}
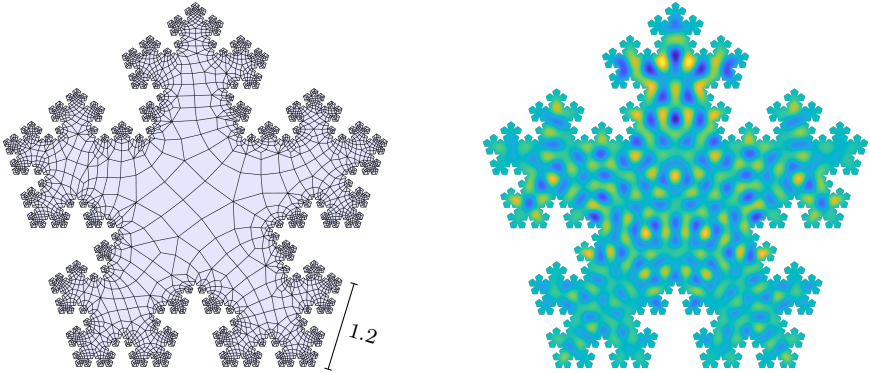

As a practical example of $hp$-adaptivity, we consider using \ultraSEM on a domain with small-scale geometric features along its boundary. The domain $\Omega$ is a snowflake shape created by a fractal-like Penrose tiling (see \cref{fig:snowflake} (left)). We construct a mesh of 4,568 quadrilaterals over $\Omega$ using the meshing software Gmsh~\cite{Gmsh}, with the element size constrained to be smaller near the boundary and larger in the interior. To specify a $p$-adaptive discretization, we define a function that varies smoothly from $p=40$ in the center of $\Omega$ to $p=7$ near the boundary, indicating that coarse elements in the interior of $\Omega$ employ a high-$p$ discretization while fine elements close to the boundary of $\Omega$ employ a lower $p$. The total number of degrees of freedom for this $hp$-mesh is $N = \text{333,627}$. We locate the domain $\Omega$ such that $y<0$ for all $(x,y) \in \Omega$, and solve the gravity Helmholtz equation
\begin{equation}\label{eq:gravity_helmholtz}
\nabla^2 u + 100(1-y)u = -1, \qquad u \in \Omega,
\end{equation}
with zero Dirichlet boundary conditions. The computation in \ultraSEM takes about 75 seconds (43 seconds in the local build stage, 31 seconds in the global build stage, and 1 second in the solve stage) and consumes approximately 10GB of memory. The computed solution is shown in \cref{fig:snowflake} (right), though is only accurate to about one digit when compared in relative infinity norm to an over-resolved solution. To obtain further accuracy, the use of impedance-to-impedance maps may be necessary to avoid artificial resonances when merging operators~\cite{Gillman_15_01}. The $h$-adaptive nature of the discretization allows for the small-scale geometry of the domain boundary to be resolved without using a prohibitive number of elements, while the $p$-adaptive nature of the discretization allows for the high-degree approximation of smooth functions on coarse elements.

\subsection{Implicit time-stepping for parabolic problems}

The ability to reuse precomputed solution operators allows for efficient implicit time-stepping for parabolic problems. To demonstrate, we consider solving the variable-coefficient convection-diffusion equation on the domain $\Omega = [0,10] \times [-1,1]$ over the time span $[0,T]$,
\begin{equation}\label{eq:convdiff}
\frac{\partial u}{\partial t} = \kappa \nabla^2 u - \nabla \cdot \left(\vec{b}(x,y) u\right), \qquad u \in \Omega \times [0,T],
\end{equation}
for $T>0$, with initial condition $u(x,y,0) = e^{-4(x-1)^2-4y^2}$ and zero Dirichlet boundary conditions. This equation models the transport of a contaminant concentration in a flow. We define the diffusivity $\kappa = 0.01$ and the convective velocity $\vec{b}(x,y) = (1 - e^{\gamma x} \cos 2\pi y, \tfrac{\gamma}{2\pi} e^{\gamma x} \sin 2\pi y)$. Here, the velocity field $\vec{b}(x,y)$ is the analytical solution to the Kovasznay flow~\cite{Kovasznay_48_01}, where $\gamma = \mathrm{Re}/2 - \sqrt{\smash{\mathrm{Re}}^2/4 - 4\pi^2}$ and $\mathrm{Re}=100$ is the Reynolds number.

\begin{figure}[htb]
	\centering
    \begin{minipage}{0.49\textwidth}
    \centering
    \includegraphics[width=\textwidth]{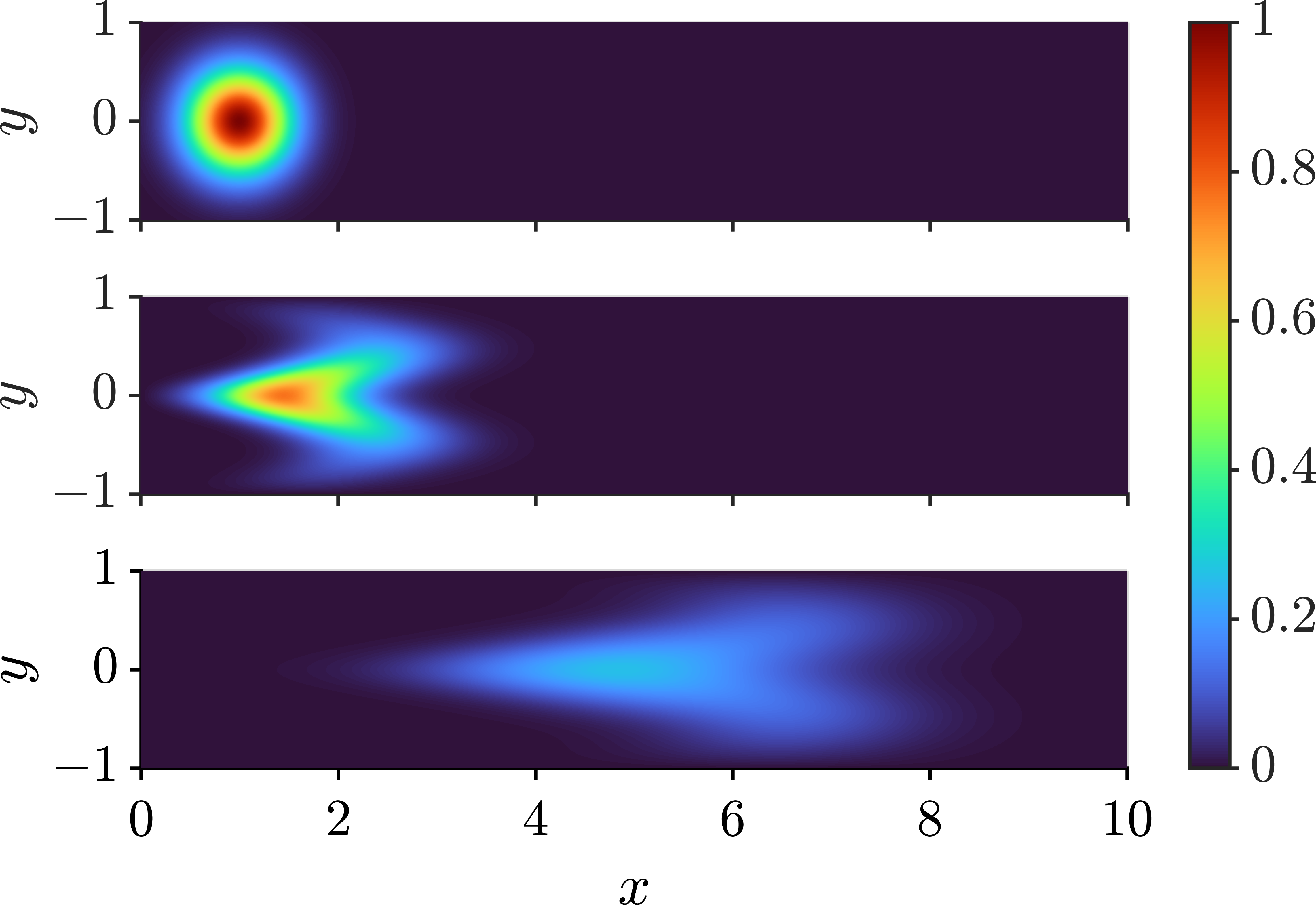}
    \end{minipage}%
    ~
    \begin{minipage}{0.49\textwidth}
    \centering
    \begin{overpic}[width=0.85\textwidth]{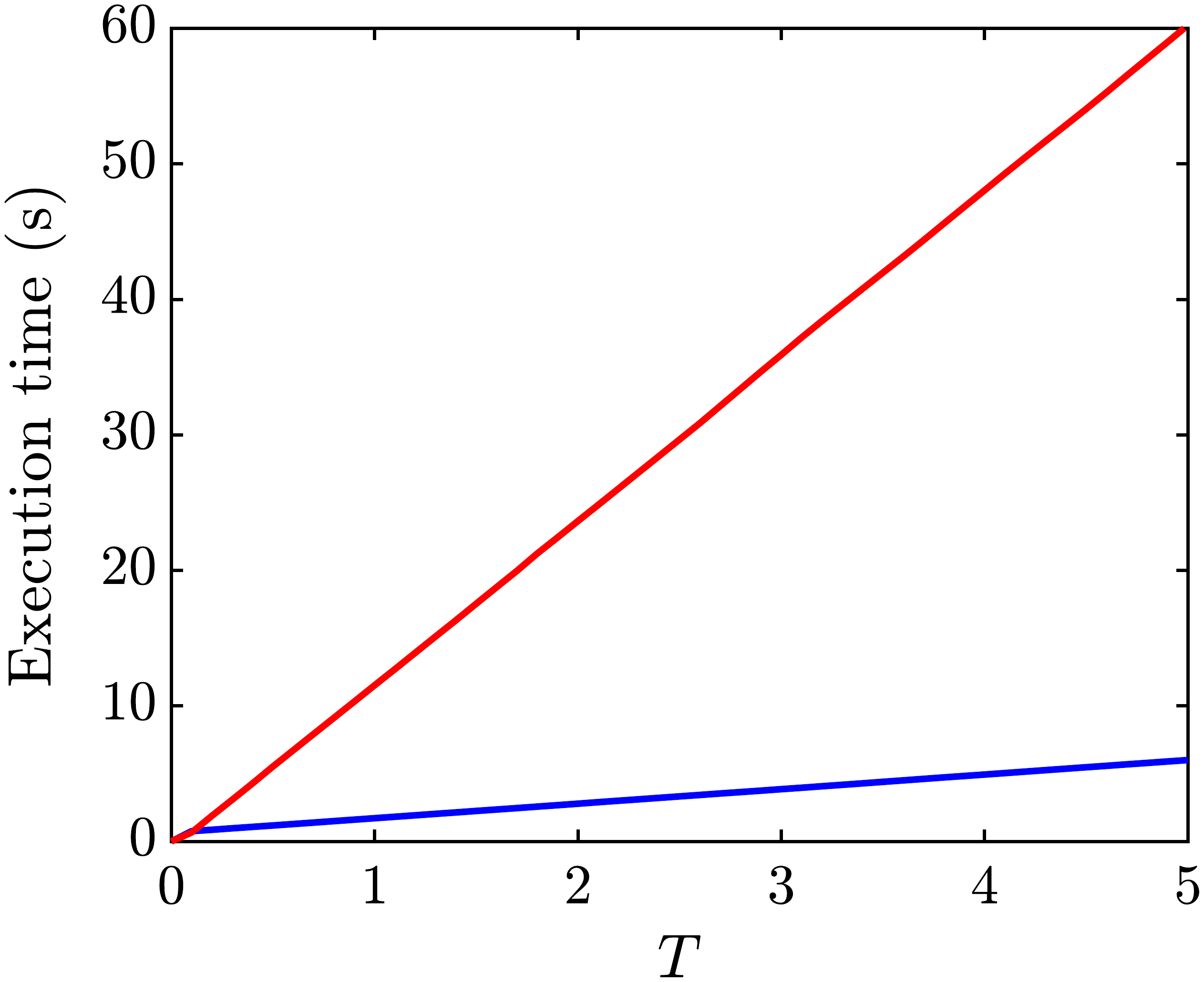}%
      \put(45,16.5) {\rotatebox{4}{\scalebox{0.65}{with \mytt{updateRHS}}}}
      \put(39,34.5) {\rotatebox{39.5}{\scalebox{0.65}{without \mytt{updateRHS}}}}
    \end{overpic}
    \end{minipage}%
\caption{(Left) Snapshots of the solution to the convection-diffusion equation \cref{eq:convdiff} at times $t=0$, $t=1$, and $t=5$, computed using \ultraSEM in space and backward Euler in time. (Right) The execution time required to simulate \cref{eq:convdiff} over the time span $[0,T]$, by either reconstructing the operators from scratch at each time step (red) or updating the particular solution using \texttt{\upshape updateRHS} (blue).}
\label{fig:timestepping}
\end{figure}

Define the time step $\Delta t = 0.1$ and time points $t_n = n \Delta t$ for integers $n \geq 0$, and let $u^n$ denote the approximate solution to \cref{eq:convdiff} at time $t_n$. Discretizing in time using the backward Euler method yields a steady-state PDE in $u^{n+1}$,
\begin{equation}\label{eq:backward_euler}
u^{n+1} - \Delta t \; \kappa \nabla^2 u^{n+1} + \Delta t \; \nabla\cdot\left(\vec{b}u^{n+1} \right) = u^n,
\end{equation}
which must be solved once per time step to compute $u^{n+1}$ from $u^n$. We use \ultraSEM to solve \cref{eq:backward_euler} on a $4 \times 20$ Cartesian mesh of $\Omega$ with polynomial order $p=16$ on each element, which yields an infinity-norm relative error of $10^{-6}$ at time $t=5$ when compared to an over-resolved solution. \cref{fig:timestepping} (left) shows snapshots of the computed solution at times $t=0$, $t=1$, and $t=5$. As the righthand side of \cref{eq:backward_euler} depends on $n$, the operators in \ultraSEM must be updated at each time step. If the operators are reconstructed from scratch at each time step, simulating to time $T=5$ completes in roughly one minute (see \cref{fig:timestepping} (right, red)). If instead only the particular solution is reconstructed using \mytt{updateRHS}, then the same simulation completes in roughly 6 seconds (see \cref{fig:timestepping} (right, blue)). \cref{fig:timestepping} (right) compares the execution times required to simulate \cref{eq:convdiff} over the time span $[0,T]$ using these two methods. It is clear that when many time steps are taken, \mytt{updateRHS} should always be used.

\section{Future work}
The development of fast direct solvers for three-dimensional problems is an active area of research~\cite{Hao_16_01}, and it may be possible to generalize the ultraspherical spectral element method to three-dimensional meshed geometries. The collocation-based HPS scheme in three dimensions has a computational complexity of $\mathcal{O}(p^9)$ due to the inversion of dense $p^3 \times p^3$ matrices on each element in the local build stage. We believe that using the ultraspherical spectral method on each leaf would reduce this complexity to $\mathcal{O}(p^7)$ or even $\mathcal{O}(p^6)$ for certain problems. However, a careful analysis of the storage costs in three dimensions is necessary to determine if such a direct solver is practical.

\section*{Acknowledgements}
We would like to thank Sheehan Olver, Keaton Burns, and Marc Gilles for their useful discussions on domain decomposition with ultraspherical polynomials, and Federico Fuentes for his expertise on $hp$-adaptivity theory. We are grateful to Heather Wilber, Sheehan Olver, Patrick Farrell, and Alex Barnett for their comments on a draft of this work.

\bibliographystyle{siamplain}
\bibliography{references}

\end{document}